\scrollmode
\documentclass[12pt]{amsart}
\usepackage{verbatim}
\usepackage{cite}
\usepackage{upref}
\usepackage{xspace}
\usepackage{eucal}
\usepackage{amssymb}
\usepackage{graphicx}
\usepackage{subfigure}
\usepackage{psfrag}
\usepackage{rcs}
\usepackage{mathrsfs}
\usepackage[all]{xy}
\SelectTips{cm}{}
\newif \ifwide
\newif \ifavnermargin
\def \makemargins{
\ifwide
	\oddsidemargin .25in
	\evensidemargin .25in
	\textwidth 6.00in
\else
\fi
\ifavnermargin
	\headheight=7pt
	\textheight=574pt
	\textwidth=432pt
	\topmargin=14pt
	\oddsidemargin=18pt
	\evensidemargin=18pt
\else	
\fi
}

\avnermargintrue
\makemargins
\theoremstyle{plain}

\theoremstyle{definition}

\theoremstyle{remark}

\RCS $Revision: 1.12 $

\newcount\timehh\newcount\timemm
\timehh=\time 
\divide\timehh by 60 \timemm=\time
\newcommand{\draftauthor}[1]{\author{#1
    {
      --- \protect \protect\sc\today\ ---
      \ifnum\timehh<10 0\fi\number\timehh\,:\,\ifnum\timemm<10 0\fi\number\timemm
      \protect \, \, \protect \bf v. \RCSRevision
    }
  }
}
\newcommand{\R}{{\mathbb R}}

\newcommand{\C}{{\mathbb C}}
\newcommand{\Z}{{\mathbb Z}}

\newcommand{\Q}{{\mathbb Q}}

\newcommand{\V}{{\mathscr V}}
\newcommand{\sB}{{\mathscr B}}
\newcommand{\Proj}{{\mathbb P}}

\renewcommand{\O}{{\mathscr O}}
\newcommand{\PPP}{\mathscr P}
\newcommand{\MMM}{\mathscr M}

\newcommand{\sC}{\mathscr C}

\newcommand{\bD}{{\mathbf{D}}}

\newcommand{\bG}{{\mathbf{G}}}

\newcommand{\bL}{{\mathbf{L}}}
\newcommand{\bM}{{\mathbf{M}}}
\newcommand{\bN}{{\mathbf{N}}}

\newcommand{\bP}{{\mathbf{P}}}
\newcommand{\bQ}{{\mathbf{Q}}}
\newcommand{\bR}{{\mathbf{R}}}
\newcommand{\bS}{{\mathbf{S}}}
\newcommand{\bT}{{\mathbf{T}}}

\newcommand{\fH}{\mathfrak H}
\newcommand{\fg}{\mathfrak g}
\newcommand{\fa}{\mathfrak a}
\newcommand{\fn}{\mathfrak n}

\renewcommand{\emptyset}{\varnothing}

\newcommand{\abcd}[4]{\left(
        \begin{smallmatrix}#1&#2\\#3&#4\end{smallmatrix}\right)}

\DeclareMathOperator{\SL}{SL}
\DeclareMathOperator{\GL}{GL}
\DeclareMathOperator{\SO}{SO}

\DeclareMathOperator{\U}{U}
\DeclareMathOperator{\SU}{SU}
\DeclareMathOperator{\Sp}{Sp}

\newcommand{\Vor}{Vorono\v{\i}\xspace}

\newcommand{\cd}[4]{\xymatrix{#1\ar[d]\ar[r]&#2\ar[d]\\ #3\ar[r]&#4}}
\newcommand{\Neq}{\stackrel{N}{\sim}}
\newcommand{\Leq}{\stackrel{L}{\sim}}
\newcommand{\Req}{\stackrel{R}{\sim}}
\newcommand{\RLeq}{\stackrel{RL}{\sim}}

\newcommand{\ol}[1]{\overline{#1}}

\newcommand{\borelserre}{\text{BS}}
\newcommand{\reductiveborelserre}{\text{RBS}}
\newcommand{\satake}{\text{Sat}}
\newcommand{\bailyborel}{\text{BB}}
\newcommand{\tits}{\text{T}}
\newcommand{\geo}{\text{geo}}
\newcommand{\martin}{\text{Mar}}
\newcommand{\Rgeqz}{\R_{\geq 0}}
\newcommand{\scal}[2]{\langle #1, #2\rangle}

\DeclareMathOperator{\Span}{span}
\DeclareMathOperator{\Sym}{Sym}

\DeclareMathOperator{\Eis}{Eis}

\DeclareMathOperator{\Int}{Int}

\CompileMatrices

\numberwithin{equation}{subsection}
\numberwithin{figure}{section}
\numberwithin{table}{section}

\begin{document}

\newif \ifdraft
\def \makeauthor{
\ifdraft
	\draftauthor{Paul E. Gunnells}
\title{Robert MacPherson and arithmetic groups (\RCSRevision)}
\else
\title{Robert MacPherson and arithmetic groups}
\author{Paul E. Gunnells}
\address{Department of Mathematics and Statistics\\
University of Massachusetts\\
Amherst, MA 01003\\
USA}
\email{gunnells@math.umass.edu}
\dedicatory{Dedicated to Robert MacPherson on the occasion of his sixtieth birthday}
\fi
}

\draftfalse
\makeauthor

\ifdraft
	\date{\today}
\else
	\date{February 19, 2006}
\fi



\begin{abstract}
We survey contributions of Robert MacPherson to the theory of
arithmetic groups.  There are two main areas we discuss: (i)
explicit reduction theory for Siegel modular threefolds, and (ii)
constructions of compactifications of locally symmetric spaces.  The
former is joint work with Mark McConnell, the latter with Lizhen Ji.  
\end{abstract}
\maketitle
\tableofcontents


\section{Introduction}\label{s:i}

\subsection{}
Arithmetic groups sit at the crossroads of many areas of mathematics.
On one hand, they lead to beautiful manifolds with intricate geometry,
and to moduli spaces for many important objects in arithmetic.  On the
other, they are conjoined with the theory of automorphic forms, and
provide one path to understanding the mysteries of the absolute Galois
group of the rationals.  The study of arithmetic groups is a
beautiful blend of algebraic topology, algebraic and differential
geometry, representation theory, and number theory, and includes some
of the most fascinating and inscrutable phenomena in mathematics.

\subsection{}
In this survey we discuss Robert MacPherson's contributions to
arithmetic groups.  We focus on two of his collaborations.  The first
(\S\ref{s:smt}) is joint work with Mark McConnell, and appears in the
papers \emph{Explicit reduction theory for Siegel modular threefolds}
(Inv. Math \textbf{111} (1993)) \cite{mmc1} and \emph{Classical
projective geometry and modular varieties} (Proceedings of JAMI 1988)
\cite{mmc2}.  The second (\S\ref{s:colss}) is joint work with Lizhen
Ji, and appears in \emph{Geometry of compactifications of locally
symmetric spaces} (Ann. Inst Fourier, Grenoble \textbf{52} (2002))
\cite{ji.macp}.  We also provide (\S\ref{s:b}) some background on
arithmetic groups and their relationship to number theory and
automorphic forms.

\subsection{}
There are two other contributions of MacPherson to arithmetic groups
that we unfortunately will not discuss, the \emph{topological trace
formula} and the \emph{geometric approach to the fundamental lemma}.
To compensate for this omission, we say a few words here.

In \cite{arthur.hecke}, using his trace formula, Arthur derived an
expression for the Lefschetz number of the action of a Hecke
correspondence on the $L^2$-cohomology of a modular variety.  In view
of the Zucker conjecture (proved by Saper and Stern \cite{saper.stern}
and, independently, by Looijenga \cite{looijenga}), which identifies
this $L^2$-cohomology with the intersection cohomology of the
Baily--Borel compactification (\S\ref{ss:zoo}), Casselman and Arthur
asked whether Arthur's formula could be interpreted as a Lefschetz
fixed point formula for intersection cohomology, and in particular
whether each term in Arthur's formula might correspond to a single
fixed point component.  The ingredients (volumes of centralizers,
orbital integrals, and averaged discrete series characters) in
Arthur's formula did not, initially, look like the local contributions
one might expect from a Lefschetz fixed point formula.

Nevertheless, Arthur and Casselman's suggestion was eventually
realized in a series of papers \cite{lefschetz, announcement,
weighted, toptrace, discrete, discrete.correction} in which it was
shown that each term in Arthur's formula corresponds to a sum of
(Lefschetz) contributions over a certain collection of fixed points.
New techniques were developed including (i) a general topological
formula for the local contribution to the Lefschetz number from
``hyperbolic" fixed points; see also the closely related results of
Kashiwara and Schapira \cite{kashiwara.schapira}; (ii) the ``weighted
cohomology" of modular varieties (similar to the intersection
cohomology, but involving a weight truncation rather than a degree
truncation), later shown \cite{nair} to coincide with the
\emph{weighted $L^2$-cohomology} of Franke \cite{franke}; and (iii) a
combinatorial formula for the characters of discrete series, somewhat
different from the formulas \cite{herb} of Herb.  For more details
about the actual constructions involved in this program, we refer to
the announcement \cite{announcement}.

For the geometric approach to the fundamental lemma, we have
considerably less to say.  This program is still under development,
with some publications available \cite{gkm2, gkm3} and others in
preparation.  Recently Laumon--Ngo extended and improved the original
strategy of Goresky--Kottwitz--MacPherson and completed a proof of the
fundamental lemma for unitary groups, a spectacular advance
\cite{laumon.ngo}.  It is premature to predict where all this will
lead, but some experts share a certain optimism.  In any case, writing
a detailed survey of this subject would require a long and complicated
article of its own.

\subsection{Acknowledgments}
I thank Avner Ash, Mark Goresky, Lizhen Ji, and Mark McConnell for
helpful comments.  I also thank the editors for inviting me to write
this article.  Finally, I'm glad to thank (in print) Bob MacPherson,
not only for helping to create such beautiful mathematics, but also
for being a patient and benevolent teacher.


\section{Arithmetic groups}\label{s:b}

\subsection{}
In this section we review background on algebraic groups, arithmetic
groups, and discuss their relationships to geometry and number theory.

\subsection{Linear algebraic groups}
\subsubsection{}
The starting point is $\bG$, a connected linear algebraic group
defined over $\Q$.  For our purposes, this means the following
\cite[\S2.1.1]{plat.rapin}:
\begin{enumerate}
\item The group $\bG$ has the structure of an affine algebraic variety
given by an ideal $I$ in the ring $\C [x_{ij}, D^{-1} ]$, where
the variables $\{ x_{ij} \mid 1\leq i,j\leq n\}$ should be interpreted
as the entries of an ``indeterminate matrix,'' and $D$ is the
polynomial $\det (x_{ij})$.  Both the group multiplication $\bG \times
\bG \rightarrow \bG$ and inversion $\bG \rightarrow \bG$ are required
to be morphisms of algebraic varieties.
                                                                                                                              
The ring $\C [x_{ij}, D^{-1} ]$ is the coordinate ring of the
algebraic group $\GL_{n}$.  Hence we are viewing $\bG$ as a subgroup
of $\GL_{n}$ defined by the vanishing of polynomial equations in matrix
entries.
\item \emph{Defined over $\Q$} means that the ideal $I$ is generated
by polynomials with rational coefficients.
\item \emph{Connected} means that $\bG $ is connected as an algebraic
variety.
\end{enumerate}

Given any subring $R \subset \C$, we can consider the group of
$R$-rational points $\bG (R)$.  As a set $\bG (R)$ is the set of
solutions to the defining equations for $\bG$ with coordinates in $R$.
Especially important in the following will be the groups of \emph{real
points} $\bG (\R)$, \emph{rational points} $\bG (\Q)$, and
\emph{integral points} $\bG (\Z)$.  Following a usual convention, we
denote algebraic groups by bold letters and their groups of real
points by the corresponding Roman letter.  Hence we write $G = \bG
(\R)$, $P = \bP (\R)$, and so on.

\subsubsection{}
A linear algebraic group $\bG$ is \emph{reductive} if its maximal
connected unipotent normal subgroup is trivial, and \emph{semisimple}
if its maximal connected solvable normal subgroup is trivial.

The prototypical example of a reductive group is $\bG = \GL_{n}$, the
\emph{split general linear group}.  For any ring $R\subset \C$, the
group $\GL_{n} (R)$ is the group of $n\times n$ invertible matrices
with entries in $R$.

There are two basic examples of semisimple groups that will be
important for us.  The first is the \emph{split special linear group}
$\bG = \SL_{n}$.  For any ring $R\subset \C$, we have $\bG (R) =
\SL_{n} (R)$, which is the group of all $n\times n$ matrices with
entries in $R$ and with determinant equal to one.  The second is the
\emph{split symplectic group} $\Sp_{2n}$.  The group of real points
$\Sp_{2n} (\R)$ is the automorphism group of a fixed nondegenerate
alternating bilinear form on a real vector space of dimension $2n$.

For some nonsplit examples, let $F$ be an algebraic number field.
Then there is a $\Q$-group $\bG_{F}$ satisfying $\bG (\Q) = \SL_{n}
(F)$.  The group $\bG_{F}$ is constructed via the technique of
\emph{restriction of scalars} from $F$ to $\Q$; the notation is
$\bG_{F} = \bR_{F/\Q} \SL_{n}$ \cite[\S2.1.2]{plat.rapin}.  For
example, if $F$ is totally real, the group $\bR_{F/\Q }\SL_{2}$ plays
an important role in the study of \emph{Hilbert modular forms}.

\subsubsection{}
Another family of examples is provided by \emph{tori}.  A torus is a
linear algebraic group $\bT$ such that $\bT\simeq \bD_{n}$, where
$\bD_{n}\subset \GL_{n}$ is the subgroup of diagonal matrices.  We
have $\bT (\C) \simeq (\C^{\times})^{n}$.  The number $n$ is called
the \emph{(absolute) rank} of $\bT$.  A torus is said to be \emph{$\Q
$-split} if $\bT (\Q ) \simeq (\Q ^{\times})^{n}$, where $n$ is the
absolute rank of $\bT$, and if the isomorphism is defined over $\Q $.
The \emph{$\Q$-rank} of $\bG$ is defined to be the dimension of a
maximal $\Q$-split torus.  For example, the $\Q$-rank of $\SL_{n}$, or
more generally $\bR_{F/\Q} \SL_{n}$, is $n-1$.

\subsection{Symmetric spaces}
\subsubsection{}
Let $\bG$ be connected semisimple, and let $G = \bG (\R)$ be the group
of real points.  Then $G$ is a connected Lie group.  Let
$K\subset G$ be a maximal compact subgroup.  The quotient $X := G/K$
is called a \emph{(global) Riemannian symmetric space}; it is
diffeomorphic to a contractible smooth real manifold.

For example, if $\bG$ is the special linear group $\SL_{n}$, then $G$
is the Lie group $\SL_{n} (\R)$.  For a maximal compact subgroup $K$
we can take $\SO (n)$, the special orthogonal group.  This is the
group of $n\times n$ real matrices $A$ with $A^{-1} = A^{t}$ and $\det
A = 1$.  Thus $X = \SL_{n} (\R)/\SO (n)$.  It is easy to compute that
$\dim X = n (n+1)/2-1$.  The most familiar example from this family of
symmetric spaces is $n=2$.  We have $\SL_{2} (\R)/\SO (2)\simeq \fH$,
where $\fH\subset \C$ is the \emph{upper halfplane} of complex numbers
with positive imaginary part.

If $G = \Sp_{2n} (\R)$, then $K = \U (n)$, the \emph{unitary group} of
all $n\times n$ complex matrices $A$ with $A^{*} = A^{-1}$, where the
star denotes conjugate transpose.  One can show that the symmetric
space $X$, called the \emph{Siegel upper halfspace}, has real
dimension $n^{2}+n$.

\subsection{}
The nonsplit examples lead to more complicated geometry.  Consider
$\bG = \bR_{F/\Q }\SL_{2}$, where $F$ is a quadratic extension of
$\Q$.  Any such field has the form $F= \Q (\sqrt{d})$, where $d\not
=0,1$ is squarefree.  There are two cases to consider, $F$ real and
$F$ imaginary.

If $F$ is real, then $d>0$.  There are two ways
to embed $F$ as a subfield of $\R$, corresponding to the two choices
for the square root of $d$.  Denote these embeddings by $a\mapsto
a^{(i)}$, $i=1,2$.  We have $G = \bG (\R) = \SL_{2} (\R) \times
\SL_{2} (\R)$.  A maximal compact subgroup is $K = \SO (2)\times
\SO (2)$, and the symmetric space is $X = \fH \times \fH$.

Now suppose $F$ is imaginary.  This time $G = \SL_{2} (\C)$
and $K=\SU (2)$; the symmetric space $X$ is the hyperbolic $3$-space
$\fH_{3}$.  

\subsection{Locally symmetric spaces}\label{ss:coho}
\subsubsection{}
Let $\Gamma \subset \bG (\Q)$ be an \emph{arithmetic group}.  By
definition this is a subgroup commensurable with $\bG (\Z)$, which
means  $\Gamma \cap \bG (\Z)$ has finite index in both $\Gamma$
and $\bG (\Z)$.  For example, $\Gamma = \SL_{n} (\Z)\subset \SL_{n}
(\Q)$ is arithmetic, as is any conjugate $g\Gamma g^{-1} $, where
$g\in \SL_{n} (\Q)$.  

\subsubsection{}\label{ss:leftaction}
The left action of $\Gamma$ on the symmetric space $X$ is properly
discontinuous.  Suppose $\Gamma$ is torsion-free.  Then $\Gamma
\backslash X$ is a manifold, and is called a \emph{locally symmetric
space}. 

\subsection{Cohomology of arithmetic groups}
\subsubsection{}\label{sss:rat.repn}
We continue to assume that $\Gamma$ is torsion-free.
Since $X$ is contractible, $\Gamma \backslash X$ is an
Eilenberg--Mac~Lane space for $\Gamma $.  In particular $\pi_{1}
(\Gamma \backslash X) \simeq \Gamma$, and all other homotopy groups
vanish.  The group cohomology of $\Gamma$ (with trivial
complex coefficients) is isomorphic to the complex cohomology of the
quotient $\Gamma \backslash X$:
\begin{equation}\label{eq:gpcohodef}
H^{*} (\Gamma ; \C) \simeq H^{*} (\Gamma \backslash X; \C).
\end{equation}
In fact \eqref{eq:gpcohodef} remains true even if $\Gamma$ has
torsion, since we are using complex coefficients.(\footnote{The
isomorphism \eqref{eq:gpcohodef} is not true for \emph{integral}
coefficients if $\Gamma$ has torsion.  Similar remarks apply to more
general coefficient modules $M$ than those we consider in this paper.})

More generally, let $M$ be a complex finite dimensional rational
representation of $G$.  Then $M$ is also a $\Gamma$-module, and one
can define the group cohomology $H^{*} (\Gamma ; M)$ of $\Gamma$ with
coefficients in $M$.  One can associate to $M$ a locally constant
sheaf $\MMM$ on $\Gamma \backslash X$ so that
\[
H^{*} (\Gamma ; M) \simeq H^{*} (\Gamma \backslash X; \MMM).
\]
Again, this remains true even if $\Gamma$ has torsion, although we
must use locally constant sheaves on orbifolds.

\subsubsection{}\label{ss:modularcurve}
It turns out that the cohomology groups $H^{*} (\Gamma ; M)$ are
intimately related to automorphic forms.  The simplest interesting
example is $\bG = \SL_{2}$.  The basic arithmetic group here is $\bG
(\Z) = \SL_2 (\Z )$, which acts on the symmetric space $\fH$ by
fractional linear transformations:
\[
\left(\begin{array}{cc}
a&b\\
c&d
\end{array} \right)\cdot z  = \frac{az+b}{cz+d}, \qquad \left(\begin{array}{cc}a&b\\
c&d
\end{array} \right)\in \SL_{2} (\Z), \quad z \in \fH.
\]

The group $\SL_{2} (\Z)$ is the most obvious arithmetic group in
$\SL_{2} (\Q)$, but there are many others.  For any $N\geq 1$, let
$\Gamma_{0} (N)\subset \SL_{2} (\Z)$ be the subgroup of matrices that
are upper triangular modulo $N$.  The locally symmetric space $Y_{0}
(N) = \Gamma_{0} (N)\backslash \fH $ is the \emph{(open) modular curve
of level $N$.}  Topologically, $Y_{0} (N)$ is a punctured surface with
genus roughly $N/12$.  The Eichler--Shimura isomorphism connects the
cohomology of $Y_{0} (N)$ with holomorphic modular forms.  In
particular, we have
\[
H^{1} (\Gamma_{0} (N); \C) = H^{1} (Y_{0} (N); \C) \simeq S_{2} (N)\oplus \overline{S}_{2}
(N)\oplus \Eis_{2} (N),
\]
where $S_{2} (N)$ is the space of weight 2 cuspidal modular forms of
level $N$, the bar denotes complex conjugation, and $\Eis_{2} (N)$ is
the space of weight 2 Eisenstein series of level $N$.  More generally,
if $k\geq 0$, let $M_{k} = \Sym^{k} (\C^{2})$ be the $k$th symmetric
power of the standard representation of $\Gamma_{0} (N)$, and let
$\MMM_{k}$ be the associated locally constant sheaf on $Y_{0} (N)$
Then we have, in obvious notation,
\begin{equation}\label{eq:es}
H^{1} (\Gamma_{0} (N); M_{k}) = H^{1} (Y_{0} (N); \MMM_{k}) \simeq
S_{k+2} (N)\oplus \overline{S}_{k+2} (N)\oplus \Eis_{k+2} (N) \quad k\geq 0.
\end{equation}
Hence the group cohomology for these modules is directly related to
modular forms of level $N$.

\subsubsection{}
Similar phenomena occur for our nonsplit examples.  Let $\bG =
\bR_{F/\Q} \SL_{2}$, where $F$ is the quadratic field $\Q (\sqrt{d})$.
Inside $F$ is a certain subring $\O_{F}$, the ring of algebraic
integers, that plays the same role for $F$ that the ring $\Z$ does for
$\Q$.  We have $\O_F = \Z [\sqrt{d}]$ unless $d=1\bmod 4$, in which
case $\O_{F} = \Z [(1+\sqrt{d})/2]$.  The basic arithmetic group here
is $\Gamma = \bG (\Z) = \SL_{2} (\O_{F})$.

Suppose $F$ is real.  Then $\abcd{a}{b}{c}{d}\in \Gamma$ acts on the symmetric
space $X = \fH \times \fH $ by fractional linear transformations,
where on
the factors we use the two embeddings $F\rightarrow \R $ to map
$\SL_{2} (\O_{F})$ into $\SL_{2} (\R)$:
\[
\left(\begin{array}{cc}
a&b\\
c&d
\end{array} \right)\cdot (z_{1},z_{2}) =
\left(\frac{a^{(1)}z_{1}+b^{(1)}}{c^{(1)}z_{1}+d^{(1)}},
\frac{a^{(2)}z_{2}+b^{(2)}}{c^{(2)}z_{2}+d^{(2)}}\right), \quad (z_{1},z_{2})\in \fH \times \fH.
\]
The quotient $\Gamma \backslash X$ is an algebraic surface, called a
\emph{Hilbert modular surface}.  Its cohomology can be described
explicitly in terms of Hilbert modular forms.

If $F$ is imaginary, the group $\Gamma$ is called a \emph{Bianchi
group}.  The quotient $\Gamma \backslash \fH_{3}$ is a 3-dimensional
hyperbolic orbifold, and its cohomology can also be computed in terms
of appropriate automorphic forms.

\subsubsection{}\label{ss:franke}
For a general arithmetic group $\Gamma$, the relationship between
cohomology and automorphic forms is captured by
the following deep theorem of Franke \cite{franke}.  

Assume the $ \Gamma$-module $M$ arises from a complex
finite-dimensional rational representation of $G$ as in
\S\ref{sss:rat.repn}.  Then Franke's result says that the cohomology
$H^{*} (\Gamma \backslash X; \MMM)$ can be systematically built out of
automorphic forms attached to $G$ and to certain subquotients of $G$.
Specifically, we have a decomposition
\begin{equation}\label{franke.decomp}
H^{*} (\Gamma ; M) = H^{*}_{\text{cusp}} (\Gamma ;
M)\oplus\bigoplus_{\{P \}}H^{*}_{\{P \}} (\Gamma ; M),
\end{equation}
where the direct sum is taken over the set of classes of
\emph{associate proper $\Q$-parabolic subgroups of $G$}
(\S\ref{ss:tits}).  The summand $H^{*}_{\text{cusp}} (\Gamma ; M)$
corresponds to the full group $G$, and is known as the \emph{cuspidal
cohomology}; this is the subspace of cohomology classes represented by
cuspidal automorphic forms.  The remaining summands constitute the
\emph{Eisenstein cohomology} of $\Gamma$.  In particular the summand
indexed by $\{P \}$ is built of Eisenstein series and their residues
attached to suitable cuspidal automorphic forms on the \emph{Levi
quotients} (\S\ref{sss:levi}) of elements of $\{P \}$; these are the
subquotients alluded to above.

For $G=\SL_{2} (\R)$, this decomposition is exactly the right side of
Eichler-Shimura isomorphism \eqref{eq:es}.  The cuspidal cohomology is
the subspace $S_{k+2} (N)\oplus \overline{S}_{k+2} (N)$, and the
Eisenstein cohomology is the subspace $\Eis_{k+2} (N)$.

For an exposition of Franke's result, as well as more information
about the cohomology of arithmetic groups, we highly recommend the
recent survey \cite{schwermer}.

\section{Reduction theory for Siegel modular threefolds}\label{s:smt}

\subsection{Reduction theory}\label{ss:quadfms}
\subsubsection{}\label{ss:reduction}
Let $\bG$ be a connected semisimple group, let $X = G/K$ be the
associated symmetric space, and let $\Gamma \subset \bG (\Q)$ be an
arithmetic group.  The goal of reduction theory is understanding
the action of $\Gamma$ on $X$.  In particular, one wants to find a
``nice'' fundamental domain for the action of $\Gamma$ on $X$.  This
should be an open set $\Omega \subset X$ such that
\begin{enumerate}
\item the union of the $\Gamma$-translates of the closure
$\overline{\Omega}$ is all of $X$, and
\item for all $\gamma \in \Gamma$ with $\gamma \not =1$, we have
$\gamma \Omega \cap \Omega = \emptyset$.
\end{enumerate}
Of course, what ``nice'' means is a matter of taste, but commonly this
is taken to mean that $\Omega$ is locally homeomorphic to a polytope,
and is of finite type in some sense. 

\subsubsection{}
The prototypical example is $\Gamma = \SL_{2} (\Z)$ and $X = \fH$.
The classical fundamental domain $\Omega$, shown in the left of Figure
\ref{red.fig}, is the set
\[
\Omega = \bigl\{z\in \fH  \bigm| |z|>1,
-1/2<\Im z <1/2 \bigr\}.
\]

This example is also the source of the name \emph{reduction theory.}
There is a close connection between $\fH$ and the space of binary
positive definite quadratic forms $Q (x,y) = Ax^{2} + Bxy + Cy^{2}$
(\S\ref{ss:posdefcone}).  The $\Gamma$-action on $\fH$ corresponds to
unimodular change of variables: the matrix $\abcd{a}{b}{c}{d}$ takes
$Q (x,y)$ to $Q (ax+cy, bx+dy)$.

This change of variables is a natural equivalence relation on
quadratic forms for the following reason.  Given a positive definite
binary quadratic form $Q$, we can form the theta series
\[
\Theta (Q) = \sum_{x,y\in \Z} q^{Q (x,y)}, \quad \text{where $q = \exp
(2\pi i z), \ z \in \fH$}.
\]
Writing $\Theta (Q) = \sum_{N} a_{N}q^{N}$, we see that the
coefficient $a_{N}$ is the number of integral solutions of the
equation $Q (x,y)=N$.  Hence $\Theta (Q)$ encodes all positive
integers that can be represented by $Q$, along with the multiplicity
of each representation.  If $Q$ and $Q'$ are related by the
$\Gamma$-action, then $\Theta (Q) = \Theta (Q')$; in particular $Q$
and $Q'$ represent the same set of integers, each with the same
multiplicity.

Now given $Q$, we can use the $\Gamma$-action to
find a form $Q'$ equivalent to $Q$ and such that the corresponding
point $z (Q') \in \fH$ lies in the domain $\overline\Omega $.
Moreover $Q'$ is uniquely determined if $z (Q') \in \Omega$.
The quadratic form $Q'$ is called \emph{reduced}, since its coefficients
$A, B, C$ are small.

\subsubsection{}\label{ss:unimodular}
\begin{figure}[htb]
\psfrag{-1}{$-1$}
\psfrag{1}{$1$}
\psfrag{0}{$0$}
\psfrag{inf}{$\infty $}
\begin{center}
\includegraphics[scale=0.4]{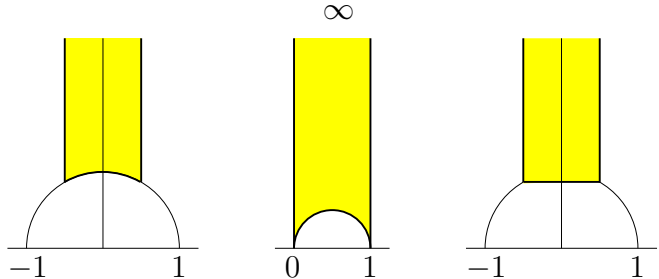}
\end{center}
\caption{Various notions of reduction domains for $\SL_{2} (\Z )$\label{red.fig}}
\end{figure}

Sometimes in reduction theory one considers weaker conditions on
$\Omega$ than those above.  For example, sometimes one replaces item
(2) in \S\ref{ss:reduction} with the requirement that $\{\gamma \mid
\gamma \Omega \cap \Omega \not =\emptyset \}$ is \emph{finite}.  The
middle of Figure \ref{red.fig} shows an example; this set is three
times as large as the exact domain on the left of Figure
\ref{red.fig}, and is stabilized by a subgroup of $\Gamma$ of order
$6$.  Although this domain is larger than a true fundamental domain,
it has the nice property that $\gamma \Omega \cap \Omega$ is
either empty or equals $\Omega$.  One may even drop this condition,
and merely require that $\{\gamma \mid \gamma \Omega \cap \Omega \not
=\emptyset \}$ is finite and that that $\Omega$ has an easily
described form.  This is useful in the theory of automorphic forms,
when one wants good control over certain analytic objects related to
$\Omega$ without fretting too much about the exact geometry of
$\Omega$.  The right of Figure \ref{red.fig} shows an example of this
kind of domain, a \emph{Siegel set} in $\fH$ \cite[I.1]{borel.ag}.

\subsection{Reduction theory and cohomology}\label{ss:gammaN}
\subsubsection{}
There is another application of reduction theory that is less
familiar:  computation of the cohomology $H^{*} (\Gamma ; M)$.
Explicit knowledge of a precise reduction domain $\Omega$ allows one
to chop $X$ up into subsets.  These subsets pass to subsets of the
quotient $\Gamma \backslash X$, and even $\Gamma '\backslash X$ where
$\Gamma '\subset \Gamma$ is a finite-index subgroup.  These subsets
can be used with standard methods of combinatorial topology to build
chain complexes to compute the cohomology of $\Gamma$ and its
finite-index subgroups.

\subsubsection{}\label{ss:princong}
We will make this discussion more precise in a moment, but the general
idea can already be seen for $\Gamma = \SL_{2} (\Z)$.  Let $N\geq 1$,
and let $\Gamma (N)\subset \SL_{2} (\Z)$ be the \emph{principal
congruence subgroup} of matrices congruent to the identity modulo
$N$.(\footnote{A \emph{congruence subgroup} of $\SL_{2} (\Z)$ is one
that contains $\Gamma (N)$ as a subgroup for some $N$.  Hence the
groups $\Gamma_{0} (N)$ from \S\ref{ss:modularcurve} are congruence
subgroups.})  For $N\geq 2$ the group $\Gamma (N)$ is torsion-free.
The quotients $Y (N) = \Gamma (N)\backslash \fH$ are punctured
topological surfaces with genus roughly $N^{3}/24$.

Let $\Omega$ be the ideal triangle in $\fH$ with vertices at
$0,1,\infty$, as in the middle of Figure \ref{red.fig}.  Then the
$\SL_{2} (\Z)$-translates of $\Omega$ pass to a ``punctured
triangulation of $Y (N)$,'' by which we mean the following. The
surface $Y (N)$ sits naturally inside a compact surface $X (N)$
obtained by filling the punctures with points, and there is a
triangulation of $X (N)$ such that the vertices are the punctures.
For example, if $N=3$ (respectively, $4$, $5$) then the genus of $Y
(N)$ is $0$, and $Y (N)$ equipped with this structure is isomorphic to
a tetrahedron (resp., octahedron, icosahedron) with punctures at its
vertices.

\subsubsection{}
\begin{figure}[htb]
\begin{center}
\includegraphics[scale=0.4, trim=0 330 0 0]{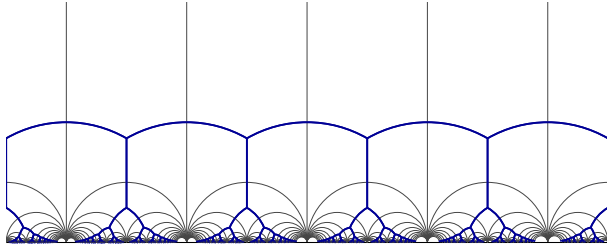}
\end{center}
\caption{The graph $W$ inside the tessellation by $\Gamma$-translates
of $\Omega $\label{tess.fig}}
\end{figure}

Here is how this structure can be used to compute $H^{*} (\Gamma (N) ;
M)$.  Consider the tessellation of $\fH$ given by the $\SL_{2}
(\Z)$-translates of $\Omega$, where $\Omega$ is as in
\S\ref{ss:princong}.  There is a regular trivalent graph $W$ embedded
in $\fH$ dual to this tessellation (Figure \ref{tess.fig}).  In other
words, $W$ has a vertex for each triangle in the tessellation, and two
vertices are connected by an edge if and only if the corresponding
triangles meet along an edge.

Modulo $\Gamma (N)$, the graph $W$ is finite.  For example, for $N=3$
(respectively, $4$, $5$), the quotient $\Gamma (N)\backslash W$ is
isomorphic to the $1$-skeleton of the tetrahedron (resp., cube,
dodecahedron).  In fact, $\Gamma (N)\backslash W$ is naturally a
deformation retract of $\Gamma (N)\backslash \fH$: the retraction is
given by pinching $Y (N)$ along its punctures to enlarge them.  Thus
the cohomology of $\Gamma (N)\backslash W$ equals that of $Y (N)$.
Since $\Gamma (N)\backslash W$ is a simplicial complex, it can be used
with standard techniques of combinatorial topology to explicitly
compute $H^{*} (\Gamma (N); M)$.  Similar considerations apply to any
finite-index subgroup $\Gamma \subset \SL_{2} (\Z)$.

\subsubsection{}
The construction of $W$ reveals more about the cohomology of $\Gamma$.
The dimension of $\fH$ is $2$, so a priori $H^{*} (\Gamma\backslash
\fH ; \MMM )$ can be nonzero in degrees $0,1,2$.  But since $\Gamma
\backslash W$ is a graph, its cohomology vanishes in degree $2$.  This
is a special case of a theorem of Borel--Serre, which implies that for
any arithmetic group $\Gamma \subset \bG (\Q)$ with $M$ as in
\S\ref{ss:franke}, we have $H^{k} (\Gamma ; M) = 0$ if $k> \nu := \dim
X - q$, where $q$ is the $\Q$-rank of $\bG$.  Thus for example if
$\Gamma \subset \SL_{n} (\Z)$, we have $\nu = n (n-1)/2$; for $\Gamma
\subset \Sp_{2n} (\Z)$ we have $\nu = n^{2}$.  The number $\nu$ is
called the \emph{virtual cohomological dimension} of $\Gamma$.  Since
for $\SL_{2} (\Z)$ the number $\nu$ is $1$, the graph $W$ is optimal
from a computational point of view: it is exactly the right dimension
to use in investigating the cohomology.  It is also optimal from an
\ae sthetic point of view, since $W$ is beautifully embedded in $\fH$.

Hence we have the following natural problem: Given an arithmetic group
$\Gamma$, find a subspace $W$ of the associated symmetric space $X$
such that the following hold:

\begin{enumerate}
\item $W$ is a locally finite regular cell complex (a regular cell
complex is a CW complex such that the attaching map from each closed
cell into the complex is a homeomorphism onto its image \cite{cooke.finney});
\item $W$ admits a cellular $\Gamma$-action such that if $\Gamma '
\subseteq \Gamma$ is a torsion-free finite-index subgroup, $\Gamma
'\backslash W$ then is a regular 
cell complex;
\item $W$ is a $\Gamma$-equivaraiant deformation retract of $X$; and 
\item $W$ has dimension $\nu (\Gamma)$.
\end{enumerate}
This is the problem
solved in \cite{mmc1, mmc2} for $\bG = \Sp_4$.

\subsection{The cone of positive-definite quadratic forms}\label{ss:posdefcone}
\subsubsection{}
To explain the results in \cite{mmc1, mmc2}, we first need to
understand how to find $W$ for the special linear group $\SL_{n}$.
Indeed, the constructions here play a key role for $\Sp_4$ through the
natural inclusion $\Sp_4 (\R )\subset \SL_{4} (\R)$.

One technique to construct $W$ for $\SL_{n}$ builds on work of \Vor's
theory of perfect quadratic forms \cite{vor1}.  We recall the
definitions.  Let $V$ be the $\R $-vector space of all symmetric
$n\times n$ matrices with entries in $\R$, and let $C\subset V$ be the
subset of positive-definite matrices.  The space $C$ can be identified
with the space of all real positive-definite quadratic forms in $n$
variables: in coordinates, if $x= (x_{1},\dotsc ,x_{n})^{t}\in \R^{n}$
(column vector), then the matrix $A\in C$ induces the quadratic form
\[
x \longmapsto x^{t}Ax.
\]
It is well known that any positive-definite quadratic form arises
in this way.  The space $C$ is a cone, in that it is preserved by
homotheties: if $x\in C$, then $\lambda x\in C$ for all $\lambda \in
\R_{>0}$.  Let $D$ be the quotient of $C$ by homotheties.  

\subsubsection{}
The case $n=2$ is illustrative.  We can take
coordinates on $V \simeq \R^{3}$ by representing any matrix in $V$ as 
\[
\left(\begin{array}{cc}
x&y\\
y&z
\end{array} \right), \quad x,y,z\in \R.  
\] 
The subset of singular matrices $Q = \{xz-y^{2}=0 \}$ is a quadric
cone in $V$ dividing the complement $V\smallsetminus Q$ into three
connected components.  The component containing the identity matrix is
the cone $C$.  The quotient $D$ can be identified with an open
$2$-disk.

\subsubsection{}
The group $G = \SL_{n} (\R)$ acts on $C$ on the left by
\[
(g,c)\longmapsto gcg^{t}.
\]
This action commutes with that of the homotheties, and thus descends
to a $G$-action on $D$.  One can show that $G$ acts transitively on
$D$ and that the stabilizer of the image of the identity matrix is
$K=\SO (n)$.  Hence we may identify $D$ with the symmetric space $X_{\SL } =
\SL_{n} (\R )/\SO (n)$.  We will do this in the sequel, using the
notation $D$ when we want to emphasize the coordinates coming from the
linear structure of $C\subset V$, and the notation $X_{\SL}$ for the
quotient $G/K$.

The $G$-action induces an action of $\SL_{n} (\Z )$ on $C$.  
This is the unimodular change of
variables action on quadratic forms as in \S\ref{ss:unimodular}.
Under our identification of $D$ with $X_{\SL}$, this is the usual action of
$\SL_{n} (\Z)$ by left translation from \S\ref{ss:leftaction}.

\subsection{The \Vor polyhedron}\label{vor.poly.section}
\subsubsection{}
Recall that a point in $\Z ^{n}$ is said to be \emph{primitive} if the
greatest common divisor of its coordinates is $1$.  In particular, a
primitive point is nonzero.  Let $\PPP\subset \Z ^{n}$ be the set of
primitive points.  Any $v\in \PPP$, written as a column vector,
determines a rank-one symmetric matrix $q (v)$ in the closure $\overline C$
via $q ( v) = vv^{t}$.  The \emph{\Vor \ polyhedron} $\Pi $ is defined
to be the closed convex hull in $\overline C$ of the points $q (v)$, as $v$
ranges over $\PPP$.  Note that by construction, $\SL_{n} (\Z )$ acts
on $\Pi $, since $\SL_{n} (\Z)$ preserves the set $\{q (v) \}$ and
acts linearly on $V$.

\subsubsection{}\label{ss:vor.fan}
The polyhedron $\Pi$ is quite complicated: it has infinitely many
faces, and is not locally finite.  However, one of \Vor 's great
insights is that $\Pi$ is actually not as complicated as it seems.
To explain his insight, we need the notion of a perfect quadratic
form.

For any $A\in C$, let $\mu (A)$ be the minimum value attained by $A$
on $\PPP$, and let $M (A)\subset \PPP$ be the set on which $A$ attains
$\mu (A)$.  Note that $\mu (A)$ is positive and $M (A)$ is finite
since $A$ is positive-definite.  Then $A$ is called \emph{perfect} if
it is recoverable from the knowledge of the pair $(\mu (A), M (A) )$.
In other words, given $(\mu (A), M (A) )$, we can write a system of
linear equations
\begin{equation}\label{system}
mZm^{t} = \mu (A), \quad m\in M (A),
\end{equation}
where $Z = (z_{ij})$ is a symmetric matrix of variables.  Then $A$ is
perfect if and only if it is the unique solution to the system
\eqref{system}. 

\subsubsection{}
We can now summarize \Vor 's main results:
\begin{enumerate}
\item There are finitely many equivalence classes of perfect forms
modulo the action of $\SL_{n} (\Z)$.  \Vor \ even gave an explicit
algorithm to determine all the perfect forms of a given dimension.
\item \label{item:minvec} The facets of $\Pi$, in other words the
codimension $1$ faces, are in bijection with the rank $n$ perfect
quadratic forms.  Under this correspondence the minimal vectors $M
(A)$ determine a facet $F_{A}$ by taking the convex hull in $\overline C$
of the finite point set $\{q (m)\mid m\in M (A) \}$.  Hence there are
finitely many faces of $\Pi$ modulo $\SL_{n} (\Z)$, and thus finitely
many modulo any finite-index subgroup $\Gamma$.
\item Let $\V$ be the set of cones over the faces of $\Pi $.  Then $\V$
is a \emph{fan}, which means (i) if $\sigma \in \V$, then any face of
$\sigma$ is also in $\V$; and (ii) if $\sigma ,\sigma '\in \V$, then
$\sigma \cap \sigma '$ is a common face of each.(\footnote{Strictly
speaking, \Vor \ actually showed that every codimension 1 cone is
contained in exactly two top dimensional cones.}) 
\item The \Vor fan $\V$ provides a reduction theory for $C$ in the
following sense: any point $x\in C$ is contained in a unique $\sigma
(x) \in \V$, and the stabilizer subgroup $\{\gamma \in \SL _{n} (\Z )
\mid \gamma \cdot \sigma (x) = \sigma (x) \}$ is finite.  \Vor\ also
gave an explicit algorithm to determine $\sigma (x)$ given $x$, the
\emph{\Vor \ reduction algorithm}.
\end{enumerate}

The number of equivalence classes of perfect forms modulo the action
of $\GL_{n}(\Z)$ grows rapidly with $n$.  \Vor computed the
equivalence classes for $n\leq 5$ \cite{vor1}.  Currently the largest
$n$ for which the number is known is $n=8$:
Dutour--Schurmann--Vallentin recently showed that there are $10916$
equivalence classes \cite{dutour}.  For a list of perfect forms up to $n = 7$, see
\cite{cs}.

\subsection{The \Vor decomposition and the retract}
\subsubsection{}
Here is how the \Vor fan $\V$ can be used to construct
higher-dimensional analogues of the tessellation in Figure
\ref{tess.fig}.  The idea is to use the cones in $\V$ to chop the
quotient $D $ into pieces.

For any $\sigma \in \V$, let $\sigma^{\circ}$ be the open cone
obtained by taking the complement in $\sigma$ of its proper faces.
Then after quotienting by homotheties, the cones $\{\sigma^{\circ }
\cap C \mid \sigma \in \V \}$ pass to locally closed subsets of $D$.
Note that $\sigma^{\circ} \cap C$ may be empty.  If $\sigma^{\circ}
\cap C \not = \emptyset$, and $c$ is the image of $\sigma^{\circ}$
modulo homotheties, we say $\sigma^{\circ}$ \emph{induces} $c$.  Note
that each $c$ is a cell, in other words is homeomorphic to an open ball,
since it is the quotient of an open polyhedral cone by homotheties.

\subsubsection{}
Let $\sC$ be the set of cells $c$ induced from the \Vor cones.
Clearly $C$ is the union of all cones $\sigma$ that induce cells in
$\sC$.  Since $\sC$ comes from the fan $\V$, the cells in $\sC$ have
good incidence properties: the closure in $D$ of any $c\in \sC$ can be
written as a finite disjoint union of elements of $\sC$.  Moreover,
$\sC$ is locally finite: by taking quotients of only the $\sigma^{\circ
}$ meeting $C$, we have eliminated the open cones lying in $\overline
C$, and it is the latter cones that are responsible for the failure of
local finiteness of $\V$.

We summarize these properties by saying that $\sC$ gives a
\emph{cellular decomposition} of $D$.  Clearly $\SL_{n} (\Z)$ acts on
$\sC$, since $\sC$ is constructed using the fan $\V$.  Thus we obtain
a cellular decomposition of $\Gamma \backslash D $ for any
finite-index $\Gamma\subset \SL_{n} (\Z)$.
We call $\sC$ the \emph{\Vor \ decomposition} of $D$.

\subsubsection{}\label{ss:buildingW}
Now we explain how $\V$ can be used to construct the cell complex $W$.
The first step is to enlarge the cone $C$ to a \emph{partial Satake
compactification} $C^{*}$.  Let $H$ be a hyperplane in $V$, and let
$\overline C$ be the closure of $C$ in $V$.  We say $H$ is a {\em
supporting hyperplane} of $C$ if $H$ is rational and $H\cap
C=\varnothing $ but $H\cap \overline C \not =\varnothing $.  Since
$\overline C$ is convex, these conditions imply that $\overline C$
lies entirely in one of the two closed half-spaces determined by $H$.
                                                                                                                              
Given a rational supporting hyperplane $H$ of $C$, let $C' =
\Int(H\cap \bar C)$, where $\Int(\phantom{a})$ denotes the interior in
the linear span.  Then $C'$ is called a {\em rational boundary
component}; it is isomorphic to a smaller dimensional cone of positive
definite quadratic forms.  Let $C^{*}$ be the union of $C$ and all
its proper rational boundary components.  We can similarly form
$D^{*}$, the union of $D$ and the images of the rational boundary
components modulo homotheties.  One can topologize $D$ such that
$\Gamma \backslash D^{*}$ is a compact Hausdorff space.  In
general $\Gamma \backslash D^{*}$ is singular.

Again, we can consider these constructions for $\SL_{2}$ and the
principal congruence subgroup $\Gamma (N)$.  In this case $D^{*}$
is $\fH \cup \Q \cup \{\infty \}$.  The quotient $\Gamma (N)\backslash
D^{*}$ is isomorphic to the surface $X (N)$ obtained by filling in
the punctures of $Y (N)$.

\subsubsection{}
It turns out that formation of the \Vor polyhedron $\Pi$ is compatible
with construction of the rational boundary components, which means the
fan $\V$ actually lies in $C^{*}$.  Let $B (\V)$ be the first
barycentric subdivision of $\V $.  Take the cones in $B (\V)$ that are
\emph{dual} to the cones inducing cells in $\sC$.  These cones are
contained entirely in $C$, and project to a collection of cells in
$D$.  By taking unions of these cells, we build $W$.  This strategy
was devised by Ash in \cite{ash77}, and used by him to construct
retracts for a large class of arithmetic groups.

\begin{figure}[htb]
\begin{center}
\includegraphics[scale=0.3]{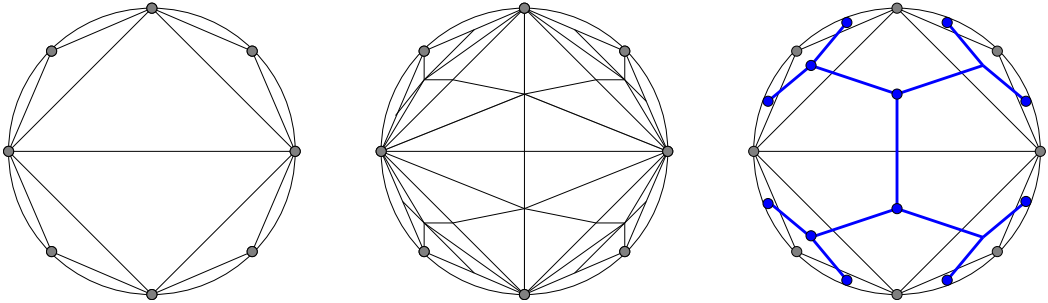}
\end{center}
\caption{Forming $W$ by taking the dual\label{buildingW.fig}}
\end{figure}

Figure \ref{buildingW.fig} shows the process in $D$ for $\SL_{2}$.  On
the left we see part of the \Vor tessellation in $D^{*}$,
displayed in the linear coordinates on $D$.  The gray discs correspond
to the points in $D^{*}\smallsetminus D$.  The middle shows the
barycentric subdivision of this decomposition.  The right shows $W$,
in heavy black lines.  Note that edges from the centers of the
triangles to the boundary components do not appear in $W$, since they
are not dual to cells in the original tessellation.  Also note that
the cells in $W$ are formed by taking unions of cells in the
barycentric subdivision.  For more about $W$, as
well as extensions of $W$ to other settings, we refer to
\cite{markmc.vor, yasaki1, yasaki2, ash.wrr, ash77, soule}.

\subsection{The symplectic cell decomposition}
\subsubsection{}
Finally we come to the symplectic case, and to the results of
\cite{mmc1, mmc2}.  Let $\bG = \Sp_{2n}$, $G =
\Sp_{2n} (\R)$, and $\Gamma = \Sp_{2n} (\Z)$.  Let $X_{\Sp}$ be the
symplectic symmetric space $\Sp_{2n}/\U (n)$, a real manifold of
dimension $n^{2}+n$.  How can we construct $W$?  Unfortunately, there
is no analogue of the cone $C$ for $X_{\Sp }$: the Siegel upper
halfspace has no hidden linear structure.

The main idea of \cite{mmc1} is extremely simple.  There is an
embedding $\Sp_{2n} (\R ) \subset \SL_{2n} (\R)$ that induces an
embedding $\iota \colon X_{\Sp} \rightarrow X_{\SL}$.  In terms of the
linear coordinates on $D \simeq X_{\SL}$, the image of $X_{\Sp}$ is
cut out by quadratic equations.  In $X_{\SL}$ we have the \Vor cells
$\sC$ that provide a cellular decomposition, and we can consider all
possible intersections $\{ c\cap \iota (X_{\Sp}) \mid c\in \sC \}$.
These provide candidates for the cells in a cellular decomposition
$\sC_{\Sp} $ of
$X_{\Sp}$. 

However, there is a potential pitfall with this idea: how do we know
that the intersections $c\cap \iota (X_{\Sp})$ are actually
\emph{cells}?  If each $c$ met the image of $X_{\Sp}$ transversely,
then we would know by a general argument that the intersections would
be cells.  Unfortunately, this is not the case.  It seems that
sometimes the \Vor cells meet the Siegel space transversely, and
sometimes not.  Thus we have no way of knowing a priori that the sets
$c\cap \iota (X_{\Sp})$ are actually nice.

\subsubsection{}
Because of this difficulty, the authors decide to study one special
case: $\Sp_{4}$.  There are several reasons why this is a good choice.
First of all, the \Vor complex is only known in full detail
for $n\leq 5$, so this is the only symplectic case other than $\Sp_2
\simeq \SL_{2}$ that can be studied explicitly.  It is also the
simplest example of a nonlinear symmetric space of $\Q$-rank $>1$.
This is significant, since to construct the symplectic retract
$W_{\Sp}$ the authors plan to use a Satake compactification as in \S
\ref{ss:buildingW}, and with such compactifications there are
significant differences in passing from $\Q$-rank 1 to $\Q$-rank $>1$.
Finally, arithmetic quotients of $\Sp_{4} (\R)$ have always played a
special role in the literature.  For one, they are moduli spaces of
abelian surfaces with extra structure, and hence give  interesting yet
tractable examples of Shimura varieties other than the modular
curves.  Also the associated automorphic forms have long been of
interest in arithmetic, since they are the first examples of Siegel
modular forms that are not elliptic.

Thus to study $\sC_{\Sp}$ the authors resort to explicit computations
in coordinates.  They find that indeed the intersections $c\cap \iota
(X_{\Sp})$ are always cells, and that the collection $\sC_{\Sp}$
provides a cell decomposition of $X_{\Sp}$.

\subsection{The symplectic retract}
\subsubsection{}
After forming the cell decomposition $\sC_{\Sp}$, McConnell and
MacPherson onstruct the symplectic retract $W_{\Sp} $ by following the
strategy indicated in \S\ref{ss:buildingW}.  They enlarge $X_{\Sp}$ to
a Satake compactification $X_{\Sp}^{*}$, and show that the
decomposition $\sC_{\Sp}^{*}$ extends to $X_{\Sp}^{*}$ in such a way
that the compact quotient $\Gamma '\backslash \sC_{\Sp}^{*}$ is a
regular cell complex for any torsion-free finite-index $\Gamma '
\subset \Sp_{4} (\Z)$.  Then they define $W_{\Sp}$ to be the
Poincar\'e dual complex to $\sC_{\Sp}$ in the first barycentric
subdivision of $\sC_{\Sp}^*$, as indicated in Figure
\ref{buildingW.fig}.

\subsubsection{}
The difficulty in carrying out these arguments is that they must first
verify that $\sC_{\Sp}^{*}$ is a regular cell complex.  Again this
doesn't follow from any general principles; they really need to show that
the closure of each cell in $\sC_{\Sp}^{*}$ is a closed ball.  Various
bad things could happen when taking the closures.  For example, the
boundaries of the closures could be non-simply-connected homology
spheres instead of true spheres.

To verify the regularity, more explicit computations are needed.  They
show that the boundaries of the closures are spheres by constructing
explicit shellings.  Roughly speaking, a shelling of a simplicial
complex is a total ordering of its maximal faces satisfying certain
properties that guarantee that the full complex is assembled from the
maximal faces nicely.  A result from combinatorics \cite{dk} states
that if a finite $n$-dimensional simplicial complex $\Delta$ is
shellable and any $(n-1)$-simplex is contained in exactly two
$n$-simplices, then $\Delta$ is homeomorphic to an $n$-sphere.  Using
a computer, McConnell and MacPherson construct shellings of the
barycentric subdivisions of the boundaries of the closures of cells
with with a computer; the largest complex they shelled was a
simplicial $S^{5}$ with $23232$ faces.

\subsubsection{}\label{sss:sp4geom}
\begin{figure}[htb]
\centering
\subfigure[Top\label{fig.crystop}]{\includegraphics[scale=0.25]{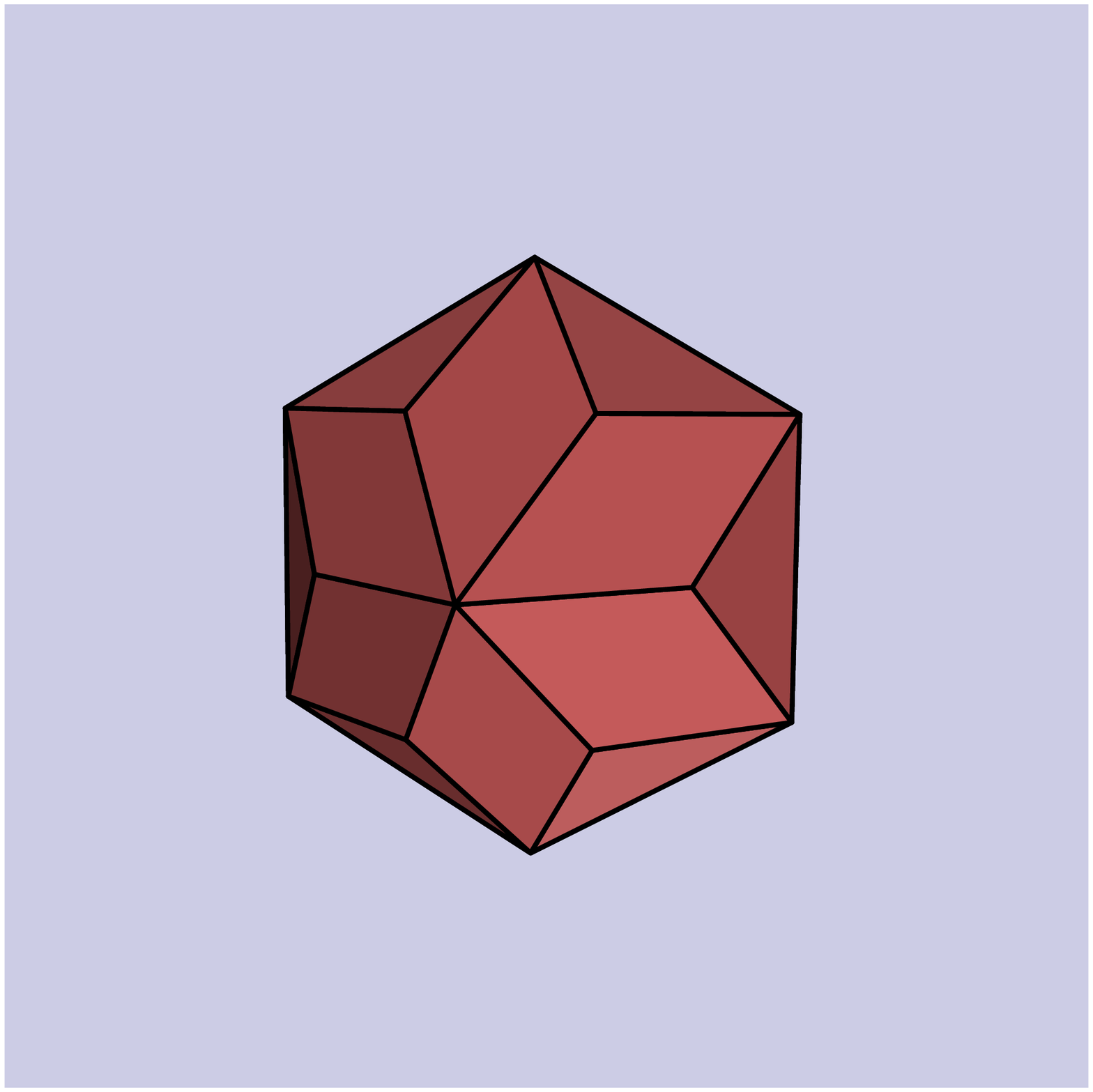}}
\quad\quad
\subfigure[Front\label{fig.crysfront}]{\includegraphics[scale=0.25]{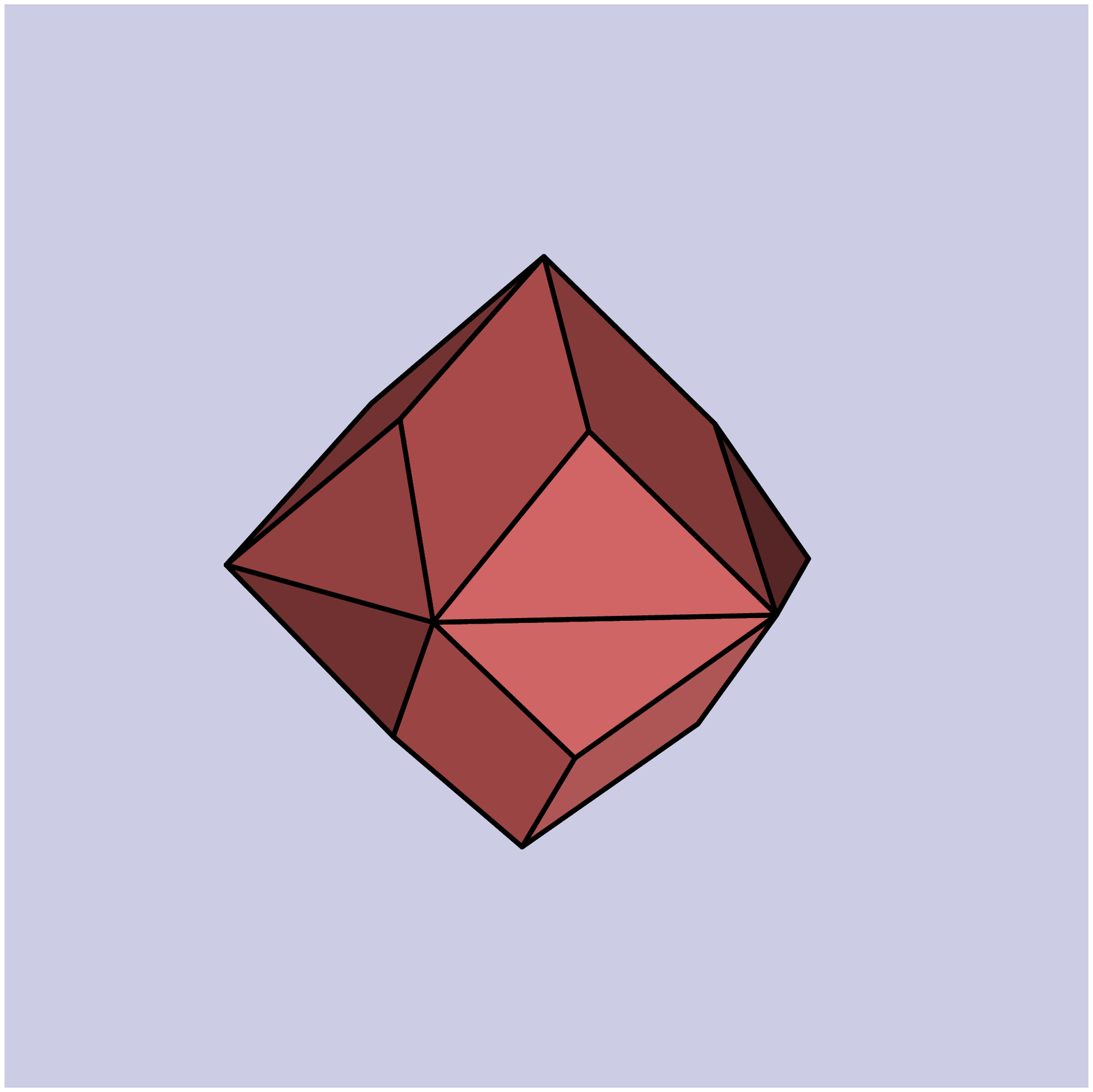}}
\caption{Two views of the crystal\label{crystal.fig}}
\end{figure}

\begin{figure}[htb]
\centering
\subfigure[Vertebra\label{fig.vert}]{\includegraphics[scale=0.25]{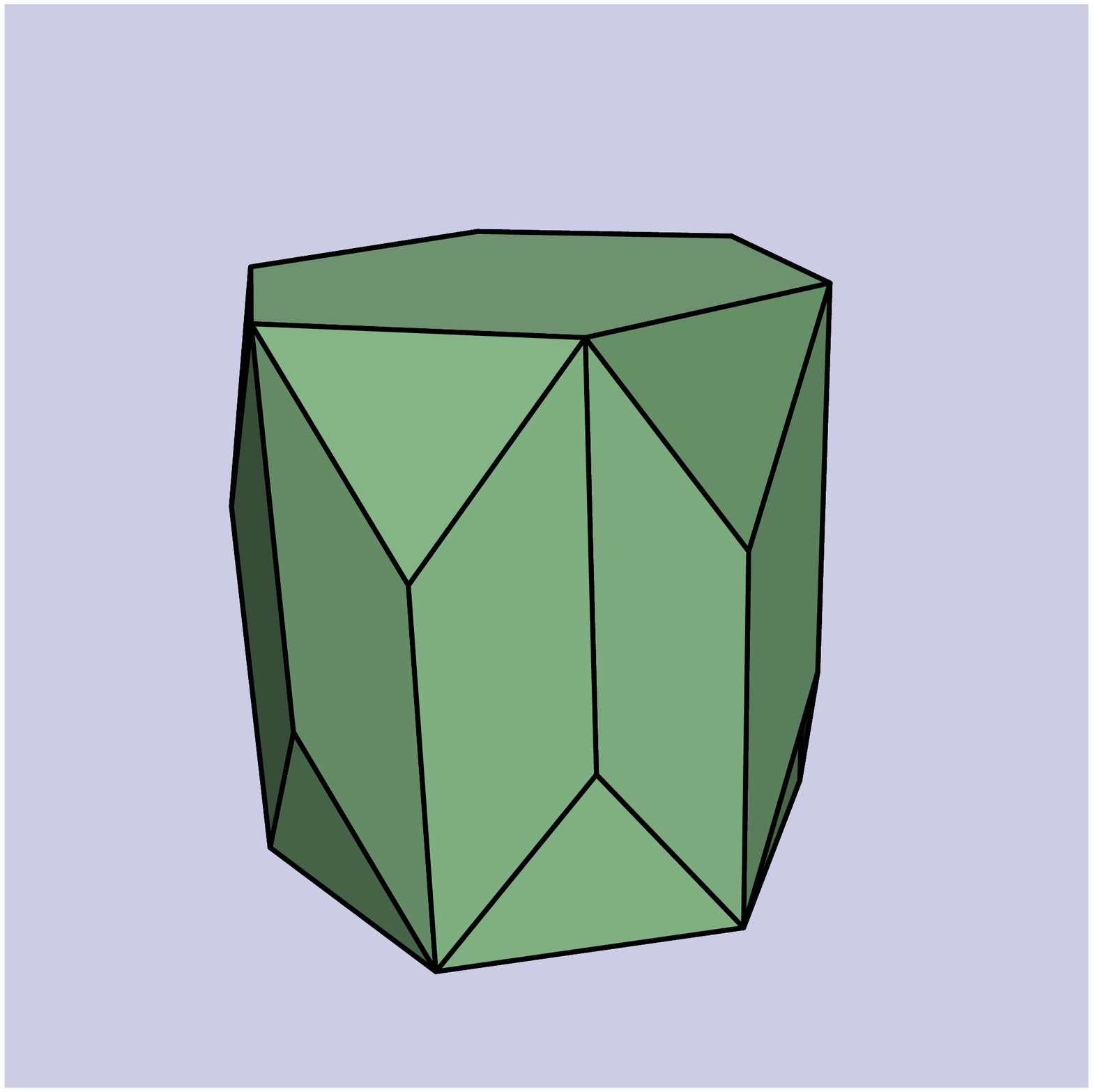}}
\quad\quad
\subfigure[Pyramid\label{fig.pyr}]{\includegraphics[scale=0.25]{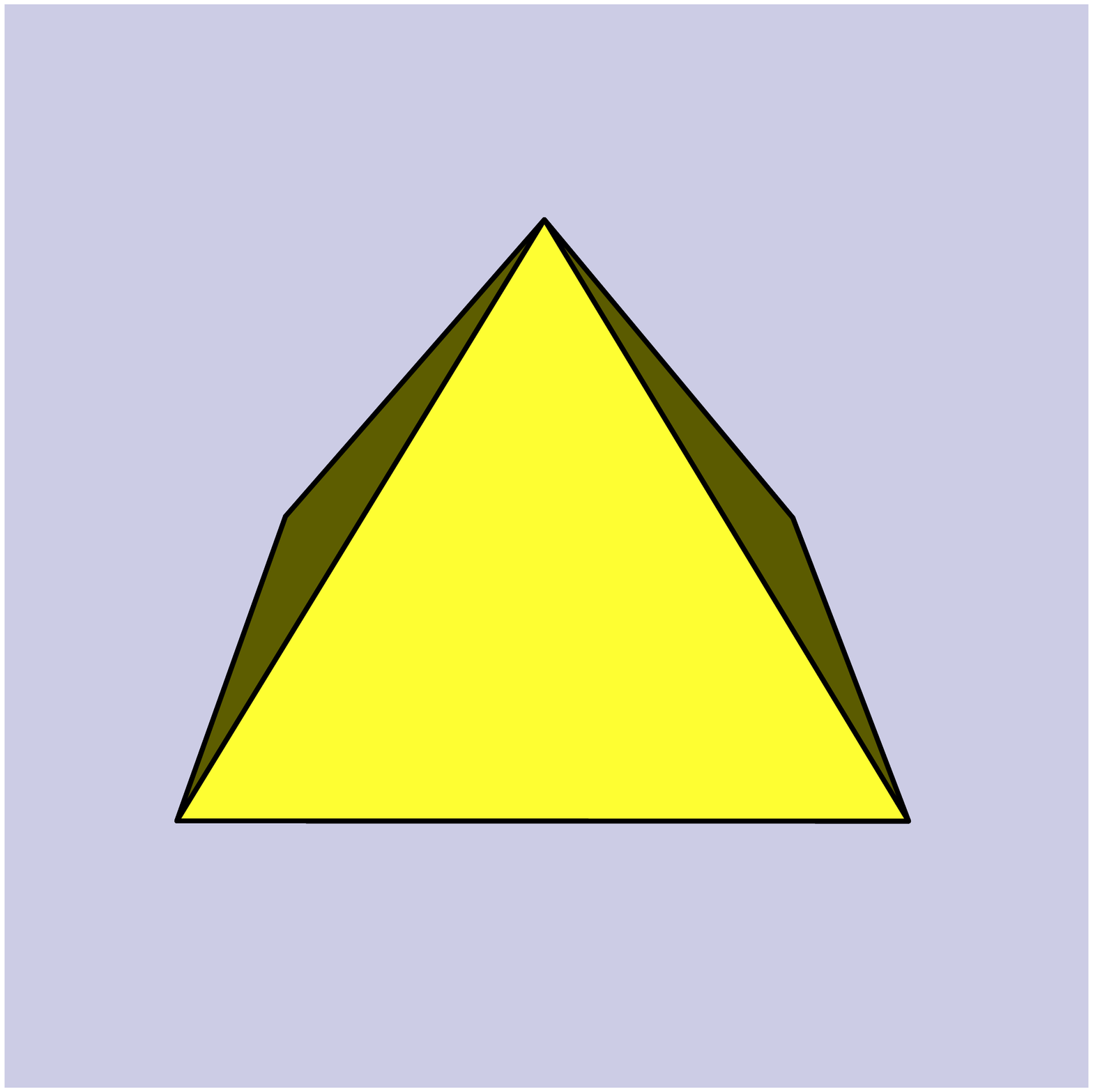}}
\caption{The vertebra and the pyramid\label{vandp.fig}}
\end{figure}

\begin{figure}[htb]
\psfrag{atinf}{at $\infty$}
\centering
\subfigure[Vertebra\label{fig.vs}]{\includegraphics[scale=0.25]{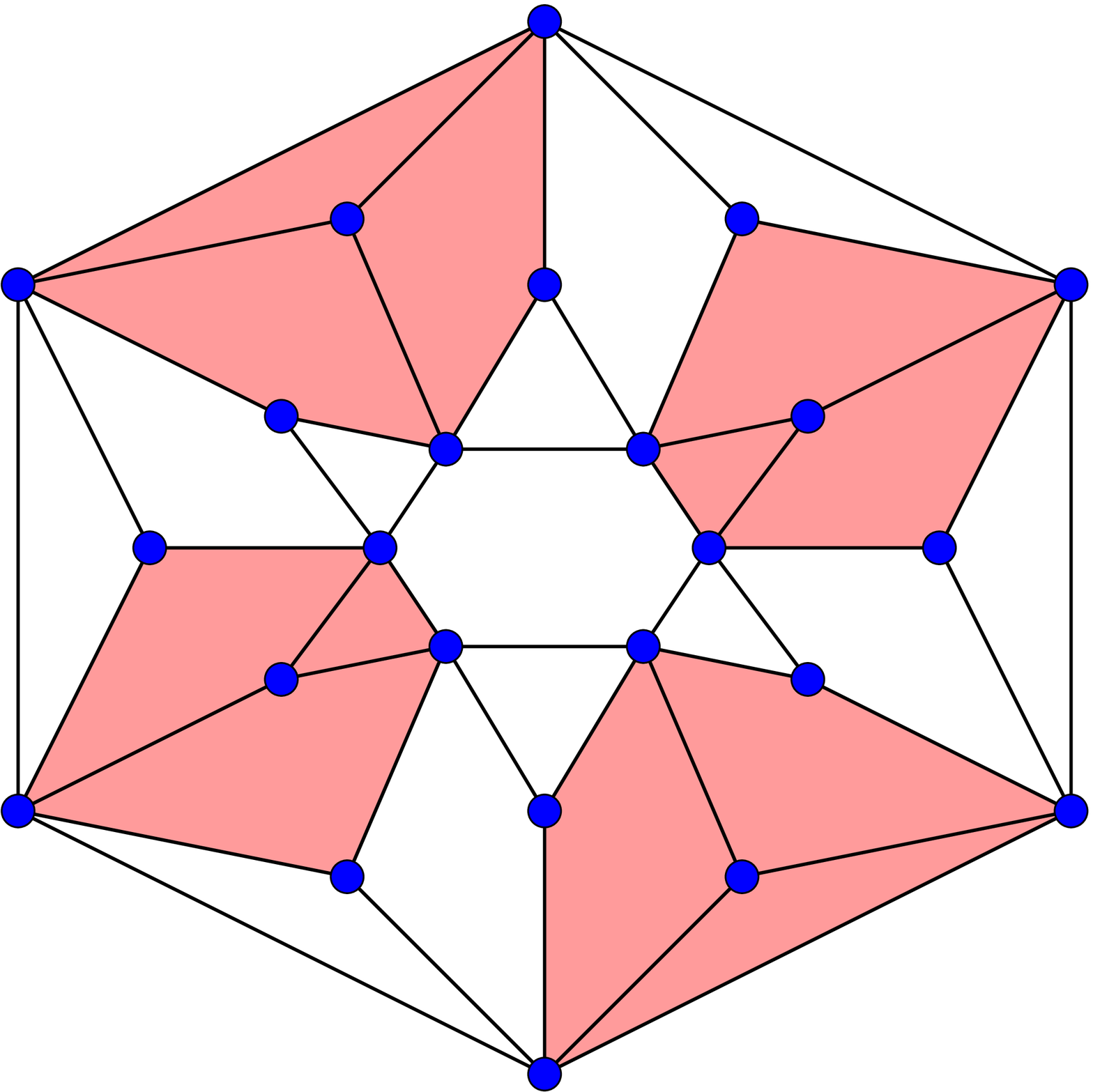}}
\quad\quad
\subfigure[Crystal\label{fig.cs}]{\includegraphics[scale=0.25]{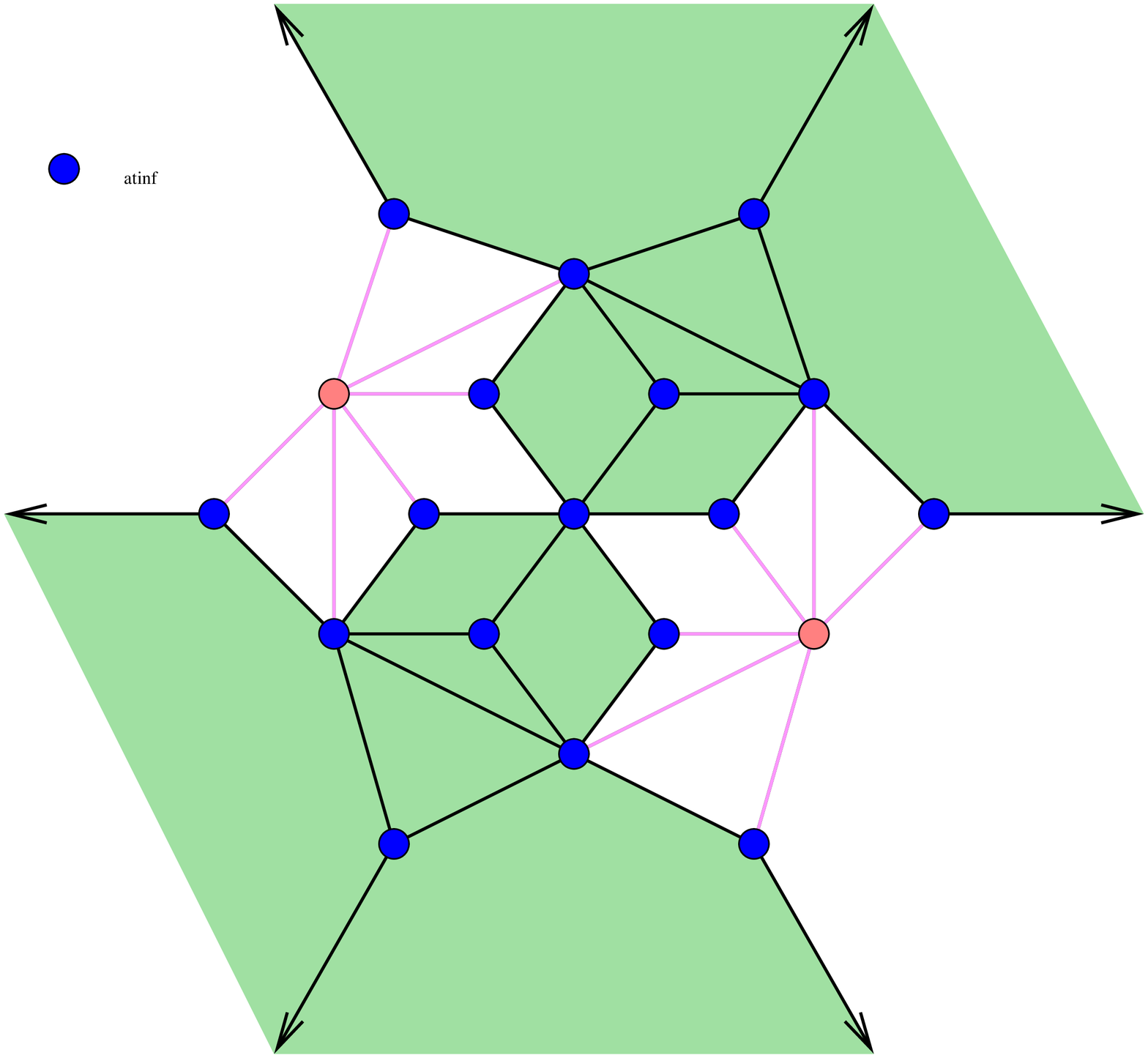}}
\caption{Schlegel diagrams for the vertebra and the crystal.  A
colored $2$-face of the vertebra (respectively crystal) also meets 
another crystal (resp., vertebra) in $\partial P$.  The light vertices
and edges in the crystal are those not contained in any vertebra
(cf. Figure \ref{fig.4cell}). 
\label{vsandcs.fig}}
\end{figure}

\begin{figure}[htb]
\begin{center}
\includegraphics[scale=0.7, trim=1500 1650 1500 800]{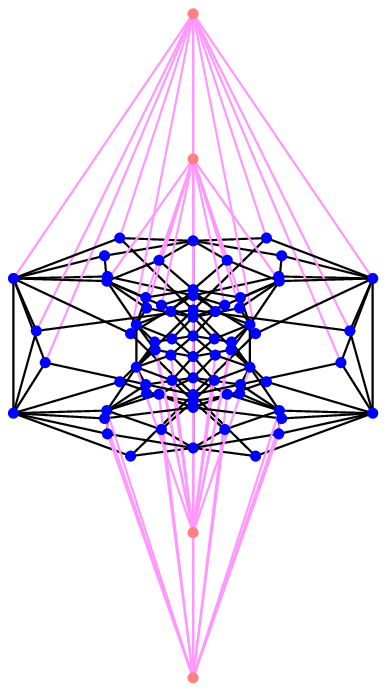}
\end{center}
\caption{The $1$-skeleton of the closure of a $4$-cell in $W_{\Sp}$.
The dark edges lie in the chain of $4$ vertebrae, which one shold
picture as being on a torus lying flat on its side. The light vertices
are the vertices of the $4$ crystals not already appearing in the
vertebrae; in $S^{3}$ they lie on a circle passing through the hole of
the torus.  If one labels these vertices $0,1,2,3$ from top to bottom,
with labels taken mod $4$, then the pair $(i,i+1)$ lies in a single
crystal (cf. the right of Figure \ref{vsandcs.fig})\label{fig.4cell}}
\end{figure}

We conclude this section by describing the geometry of the
4-dimensional cell complex $W_{\Sp}$.(\footnote{The images in Figures
\ref{crystal.fig}--\ref{fig.4cell} were produced by the author in
collaboration with Mark McConnell, as a birthday present to Bob
MacPherson.  Happy Birthday, Bob!})  The group $\Sp_{4} (\Z)$ acts
cellularly on $W_{\Sp}$, and in the following we say two cells have
the same \emph{type} if they lie in the same $\Sp_{4} (\Z)$-orbit.

There is one type of $4$-cell.  Fix a $4$-cell and let $P$ be its
closure.  Then $P$ can be realized as a cellular ball with $40$
facets, $180$ two-faces, $216$ edges, and $76$ vertices.

There are three types of $3$-cells.  In \cite{mmc2} their closures are
called the \emph{crystal}, the \emph{vertebra}, and the \emph{pyramid}
(Figures \ref{crystal.fig} and \ref{vandp.fig}).  The pyramid has $4$
triangle faces and $1$ square face.  The crystal has $12$ square and
$12$ triangle faces.  The vertebra has $2$ hexagon, $12$ square, and
$12$ triangle faces.

There are two types of $1$-cells and two types of $0$-cells modulo
$\Sp_4 (\Z)$.

Now we focus on the closure $P$ of a $4$-cell, which contains $4$
crystals, $4$ vertebrae, and $32$ pyramids.  Each crystal
(respectively, vertebra and pyramid) lies in the boundary of $3$
(resp., $3$ and $4$) $4$-cells.  Hence $P$ meets the closures of $112$
other $4$-cells.

How is $P$ assembled together from its $3$-faces?  Unfortunately $P$
has so many edges that it's difficult to draw.  The structure of $P$,
on the other hand, is relatively easy to describe.  

Begin by connecting $4$ vertebrae together in pairs along their
hexagonal faces to form a cellular solid torus.  This chain of
vertebrae has $72$ vertices, and thus accounts for all but $4$ of the
vertices of $P$.

Embed the chain of vertebrae in the boundary $\partial P \simeq S^{3}$
as a solid torus $T$, and let $C$ be an unknotted circle disjoint from
$T$ but nontrivially linking $T$ with linking number $1$.  For
example, if we identify $S^{3}$ with the polydisc $\{(x,y)\in \C^2
\mid |x| \leq 1, |y|\leq 1\}$, then we can take $T = \{|x| = 1,
|y|\leq 1 \}$ and $C = \{x=0, |y| = 1 \}$.  The remaining $4$
vertices of $P$ can be taken to lie along the circle $C$.  If we place
$4$ vertices along $C$, then any two adjacent vertices $a$, $a'$,
together with certain vertices in the vertebrae, determine a crystal.

More precisely, each vertebra meets every other crystal in $P$,
meeting each in three $2$-faces ($2$ squares and a triangle).  This is
indicated in Figure \ref{fig.vs}; the shaded $2$-faces are those that
meet the crystals.  Each triple of shaded $2$-faces meets one of
the crystals, and as we move around the chain $T$ adjacent sets of
triples meet the same crystal.  It follows that each crystal meets
every vertebra.  This is indicated in Figure \ref{fig.cs} by the
triples of shaded $2$-faces (note that in this figure, one vertex of
the crystal appears at infinity in the Schlegel diagram).

Figure \ref{fig.cs} also shows which two of the vertices of a crystal
don't lie in any vertebrae; in the figure these vertices appear
connected to the rest with light-shaded edges.  Hence to build a
crystal in $P$ one takes adjacent pairs $a$, $a'$ of vertices along
$C$ and uses each of them with $8$ vertices of the chain $T$ to flesh
out a crystal.

This shows how to find all the crystals and vertebrae in $P$.  The
pyramids simply plug up the gaps to fill out all of $\partial P$. 

Figure \ref{fig.4cell} shows the $1$-skeleton of $\partial P$.  The
torus $T$ is the central block of dark vertices and edges; the circle
$C$ is represented by a vertical line (i.e. we've used stereographic
projection to identify $S^{3}\smallsetminus \{\text{pt} \}$ with
$\R^{3}$), and the light vertices and edges are those arising from the
``fleshing out'' construction described above.

For more about the structure of $W_{\Sp}$, as well as a beautiful way
to index its cells in terms of configurations of points and lines in
$\Proj^{3}$, we refer to \cite{mmc2}.

\section{Compactifications of locally symmetric spaces}\label{s:colss}

\subsection{Introduction}
\subsubsection{}
We return to the general setting.  Let $\bG$ be a connected
semisimple $\Q$-group with group of real points $G$, and let $K$ be a
maximal compact subgroup of $G$.  Let $X=G/K$ be the associated
symmetric space, and let $\Gamma \subset \bG (\Q)$ be a torsion-free arithmetic
group.   Then as we have discussed, the quotient $\Gamma \backslash X$
is a smooth manifold.

Suppose the $\Q$-rank of $\bG$ is positive.  Then the quotient $\Gamma
\backslash X$ is noncompact.  The basic problem we now address is the
following: \emph{Can one find a good compactification of $\Gamma
\backslash X$?}

\subsubsection{}
This problem has a long history, and there are many beautiful and
ingenious solutions.  Before we talk about some of them, it's useful
say a few words about why this question is worth asking.  Why do we
care that $\Gamma \backslash X$ is noncompact? 

Indeed, for many applications we don't care.  The space $X$ carries a
natural $G$-invariant complete Riemannian metric that induces a
$G$-invariant volume form.  It turns out that even though $\Gamma
\backslash X$ is noncompact, its volume is finite (one says that
$\Gamma$ is \emph{cofinite}).  For many applications this is
sufficient; the noncompactness of $\Gamma \backslash X$ doesn't
matter. 

\subsubsection{}\label{ss:laplacian}
Nevertheless, there are reasons why one wants to compactify 
$\Gamma \backslash X$:

\begin{enumerate}
\item There are many tools from algebraic and differential topology
(for example, Lefschetz fixed point theorem, Morse theory) that only
work on compact spaces.  Since the cohomology of $\Gamma \backslash X$
is important in number theory, we would like to avail ourselves of
such tools to study it.
\item For certain groups $\bG$ the spaces $\Gamma
\backslash X$ can be interpreted as the complex points of a moduli
space, for example for $\bG  = \Sp_{2n}$.  Here the space
$\Gamma \backslash X$ parametrizes $n$-dimensional abelian varieties
with additional structure determined by $\Gamma$.  

When $\Gamma \backslash X$ is a moduli space, the points at infinity
of a compactification often correspond to
degenerations of the objects parametrized by $\Gamma \backslash X$.
The higher the $\Q$-rank, the more intricate these degenerations are.
One is interested in understanding these possible degenerations
and their interrelationships.
\item Even if $\Gamma \backslash X$ is not a moduli space, it has
additional structure coming from the group theory of $\bG$,
structure that affects computations on $\Gamma \backslash X$ in subtle
ways.  Perhaps the most definitive example of this phenomenon is
Langlands's pi\`ece de r\'esistance, the determination of the spectral
decomposition of $L^{2} (\Gamma \backslash G)$ \cite{langlands, mw}.
Other manifestations of this principle can be seen in Franke's theorem
(\S\ref{ss:franke}), Arthur's trace formula \cite{arthur.lectures},
and the computation of the spectrum of the Laplace operator on $L^{2}
(\Gamma \backslash X)$ (\S\ref{sss:laplace}).  This group-theoretic
structure also arises in the construction of compactifications of
$\Gamma \backslash X$.  Hence it is natural to understand such
computations in the fuller context of different compactifications of
arithmetic quotients of $X$.
\end{enumerate}

\subsection{Examples of compactifications}\label{ss:zoo}
\subsubsection{}
We hope the reader is convinced that compactifying
$\Gamma \backslash X$ is profitable.  Because of the wealth of
different applications of these compactifications, there is a whole
zoo of them in the literature.  Excellent overviews can be found in
\cite{borel.ji, ji.etal, goresky.lectures}.  Here we content ourselves
with a brief synopsis of some of the most useful.  To simplify
notation, we write $Y=\Gamma \backslash X$.

\subsubsection{Borel--Serre compactification $Y^{\borelserre}$ \cite{borel.serre}}
From a topological perspective, this compactification is the simplest:
it is differentiably a manifold with corners, by which we mean a
Hausdorff space smoothly modeled on the generalized halfspaces
$\R^{k}\times (\R_{\geq 0})^{n-k}$, $k=0,\dotsc ,n$.  Topologically
$Y^{\borelserre}$ is just a manifold with boundary.  The boundary
$\partial Y^{\borelserre} := Y^{\borelserre}\smallsetminus Y$ is
assembled from certain fiber bundles over locally symmetric spaces of
lower rank (\S\ref{ss:bsandrbs})

\subsubsection{Reductive Borel--Serre compactification
$Y^{\reductiveborelserre}$ \cite{zucker, weighted}}\label{sss:rbs}
In contrast to the Borel--Serre, $Y^{\reductiveborelserre}$ is
usually singular.  It is obtained as a quotient of $Y^{\borelserre}$; one
collapses the fibers in the bundles above to points.  Hence the
boundary $\partial Y^{\reductiveborelserre}$ is glued together from
locally symmetric spaces of lower rank (\S\ref{ss:bsandrbs})

Why should one ever use $Y^{\reductiveborelserre}$ instead of
$Y^{\borelserre}$?  After all, $Y^{\borelserre}$ is nearly a manifold
itself, and $Y^{\reductiveborelserre}$ is quite singular.  There are
two answers.  One is that the singularities of
$Y^{\reductiveborelserre}$ are intricate yet manageable, so little
serviceability is lost in passing from $Y^{\borelserre}$ to
$Y^{\reductiveborelserre}$.  This feature plays an important role in
the topological trace formula \cite[\S1.1.4]{toptrace}.  The other is
that $Y^{\borelserre}$ is deficient from a
differential-geometric perspective: the complete metric on $Y$ becomes
degenerate on $Y^{\borelserre}$.  This is not so on
$Y^{\reductiveborelserre}$; the metric on $Y$ extends to be
nondegenerate on $Y^{\reductiveborelserre}$.

\subsubsection{Satake compactifications $Y^{\satake}_{\tau}$
\cite{satake1, satake2}} Here there is not one compactification, but
rather a whole collection of them.  The construction is very similar
to what we did in \S\ref{ss:buildingW} in creating $C^{*}$ from $C$,
except that now we don't have a linear model for $X$.  One begins by
choosing an irreducible locally faithful representation $\tau \colon G
\rightarrow \GL (V)$ on a finite dimensional real vector space $V$.
Not just any such $\tau$ will work: $\tau$ must be \emph{geometrically
rational} (cf. \cite{casselman.rwr}).  Let $\Proj \Sym^{2}V$ be the
symmetric square of $V$ modulo homotheties.  Then one embeds the
global symmetric space $X$ in $\Proj \Sym^{2} (V)$ by $g\mapsto \tau
(g)^{t} \tau (g)$ mod homotheties.  One then takes the closure
$\overline X$ of the image of $X$ in $\Proj \Sym^{2} (V)$ and
identifies a certain subset of rational boundary components
$X^{\satake }_{\tau}\subset \ol{X}$.(\footnote{$X$ itself is also
considered to be a boundary component, called the \emph{improper}
boundary component.})  After appropriately topologizing $X^{\satake
}_{\tau }$, the quotient $Y^{\satake}_{\tau} = \Gamma \backslash
X^{\satake }_{\tau }$ is compact and Hausdorff.  This compactification
is also usually singular, in general even more singular than the reductive
Borel--Serre $Y^{\reductiveborelserre }$.

\subsubsection{Baily--Borel compactification $Y^{\bailyborel}$
\cite{baily.borel}}
This compactification, and the \emph{toroidal compactifications} that
follow, are defined for \emph{Hermitian symmetric spaces} $X$.  This
means that $X$ carries a $G$-invariant complex structure; accordingly
the quotient $Y$ does as well.  Such $X$ may be realized as a bounded
domain in a complex vector space.  One takes the closure $\overline X$
of $X$ there and again identifies a certain subset $X^{\bailyborel}$
of rational boundary components.  Again after defining an appropriate
topology, the quotient $Y^{\bailyborel} = \Gamma \backslash
X^{\bailyborel }$ is compact and Hausdorff, and in fact is a normal
analytic space.  After further investigation of the functions giving
$Y^\bailyborel$ its ringed structure, namely certain
Poincar\'e-Eisenstein series, one finds that $Y^{\bailyborel}$ is
actually a projective variety.  Topologically $Y^{\bailyborel}$ can be
obtained as one of the Satake compactifications (\cite{zucker},
cf. \cite[\S 7]{casselman.rwr}).  This compactification is usually
extremely singular.

\subsubsection{Toroidal compactifications $Y^{\Sigma}$ \cite{amrt} }
These compactifications were originally constructed to resolve the
singularities of the Baily--Borel $Y^{\bailyborel}$.  The construction
depends on some extra (noncanonical) combinatorial data $\Sigma$; this
data is a collection of rational polyhedral fans in certain
self-adjoint homogeneous cones attached to $\bG$, fans that are
geometrically very similar to the \Vor fan $\V$ in the cone of
positive definite quadratic forms $C$ (\S \ref{ss:vor.fan}).  Further
assumptions on $\Sigma$ guarantee that $Y^{\Sigma}$ is a smooth
projective variety with $\partial Y^{\Sigma}$ a divisor with normal
crossings.  An overview of the construction for Siegel modular
varieties can be found in \cite{mumford}.

\subsubsection{}
The simplest example where there are differences between some of these
compactifications is $\bG = \SL_{2}$.  Consider $Y (N)$, the modular
curve $\Gamma (N)\backslash \fH$ (\S\ref{ss:princong}).  As described
before, $Y (N)$ is homeomorphic to a punctured topological surface.
There are two obvious compactifications of $Y (N)$, one filling each
puncture with a single point, and the other lining each puncture with
a circle $S^{1}$.  The first is the reductive Borel--Serre $Y
(N)^{\reductiveborelserre}$, as well as the Baily--Borel $Y
(N)^{\bailyborel }$.  It is also the only possible Satake and only
possible toroidal compactification of $Y (N)$.  The second is the
Borel--Serre $Y (N)^{\borelserre}$.  In this case the quotient map
$\partial Y^{\borelserre}\rightarrow \partial
Y^{\reductiveborelserre}$ defined in \S\ref{sss:rbs} simply collapses
each $S^{1} \subset \partial Y^{\borelserre}$ to a point.

\subsubsection{}\label{ss:quadraticexample}
For a more revealing example, let $\bG = \bR_{F/\Q} \SL_{2}$, where
$F$ is real quadratic.  Let $\Gamma \subset \SL_{2} (\O_{F})$ be
torsion-free and of finite index.  The locally symmetric space $Y =
\Gamma \backslash \fH \times \fH$ has real dimension $4$.  Each
connected component of the boundary of the Borel--Serre $\partial
Y^{\borelserre}$ is a three manifold $Z$ that is naturally the total
space of a torus bundle over a circle: $T^{2}\hookrightarrow Z
\rightarrow S^{1}$.(\footnote{The different connected components of
$Y^{\borelserre}$ need not be homeomorphic.})

In the boundary of the reductive Borel--Serre $\partial
Y^{\reductiveborelserre}$, the tori in these bundles are collapsed to
points; hence the boundary components are circles.  This already shows
that $Y^{\reductiveborelserre}$ is singular, since the link of any
point in the boundary is a $2$-torus.(\footnote{In a stratified space
$X\subset \R^{n}$, the link $L (p)$ of a point $p$ in a stratum $S$ is
by definition $L (p) = \partial B_{\delta} (p)\cap N\cap X$, where (i)
$N$ is a submanifold through $p$ meeting all strata of $X$
transversely, and with $\dim N + \dim S = n$, and (ii) $B_{\delta}
(p)$ is a closed ball centered at $p$ with radius $0<\delta <\!\!< 1$.
It is known that the homeomorphism type of the $L (p)$ is independent
of the above choices for a wide class of stratified spaces, and the
that the normal slice $B_{\delta} (p)\cap N\cap X$ is homeomorphic to
the cone on $L (p)$.  For this same class of spaces, the link is an
invariant of the stratum $S$.  For more details, we refer to
\cite[I.1.4]{smt}.}) The Baily--Borel $Y^{\bailyborel}$ is obtained by
collapsing these circles to points; thus the link of any point in
$\partial Y^{\bailyborel}$ is a $3$-manifold $Z$ as above.  For more
details, see \cite{saper.cdmtalk}.

A toroidal compactification resolving the cusp singularities of
$Y^{\bailyborel }$ was first constructed by Hirzebruch
\cite{hirz.hms}.  Indeed, Hirzebruch's technique was one of the
principal motivations for the general theory of toroidal
compactifications \cite{amrt}.  Since $F$ is quadratic, there is 
a canonical minimal toroidal compactification.  This is not true for
higher degree $F$; compactifications in this case were first
constructed by Ehlers \cite{ehlers}.

\subsection{Building compactifications using geometry}
\subsubsection{}
The compactifications above all use the structure theory of $\bG$, the
group underlying the symmetric space $X$, in an essential way.  The
basic questions addressed by Ji--MacPherson are
\begin{itemize}
\item What natural compactifications of $Y=\Gamma \backslash X$ can be
constructed using its intrinsic geometry?  Here by intrinsic geometry
we mean objects such as sets of geodesics on $Y$.
\item What is the relationship of such compactifications with the
standard ones from \S\ref{ss:zoo}?
\end{itemize}

\subsubsection{}\label{ss:georn}
To illustrate the methods Ji--MacPherson have in mind, we
consider $M=\R^{n}$ equipped with the standard metric
\cite[\S1.1]{ji.macp}.  How can we compactify $M$?  

We first need to know what the points at infinity $M (\infty)$ should
be.  Let $m\in M$ be a fixed point, and let $\gamma \colon \Rgeqz
\rightarrow M$ be a geodesic ray with $\gamma (0)=m$.  We say $\gamma$
is \emph{based} at $m$.  Let $d (\phantom{a},\phantom{a})$ be the
distance function on $M$ induced from the metric.  Any sequence of
points $\{m_{i} \}_{i\geq 1} \subset \gamma $ with $d
(m,m_{i})\rightarrow \infty$ should converge to a point on $M
(\infty)$, since the $m_{i}$ go arbitrarily far from $m$.  Moreover,
all such sequences along $\gamma$ should tend to the same point at
infinity.

Now consider two geodesic rays $\gamma ,\gamma '$ with distinct
basepoints $m\not =m'$.  If $\gamma$ and $\gamma '$ are not parallel,
then for any two sequences $\{m_{i} \}_{i\geq 1}$ and $\{m'_{i}
\}_{i\geq 1}$ with $d (m,m_{i})$, $d (m',m'_{i})\rightarrow \infty$,
we have $d (m_{i},m_{i}')\rightarrow \infty$.  Hence these sequences
should converge to different points at infinity.  On the other hand if
$\gamma$ and $\gamma '$ are parallel, we can find two such sequences
with $d (m_{i},m'_{i})$ bounded above.  Hence it is reasonable to
require that these sequences converge to the \emph{same} point at
infinity.

This leads to the following definition.  As a set, $M (\infty)$ is the
set of geodesic rays in $M$ modulo the equivalence relation of
parallelism.  For this example the topology on $M (\infty)$ is clear:
if the angle between two rays is small, the corresponding points on
$M(\infty)$ should be close.  With this topology, $M\cup M (\infty)$
is homeomorphic to the closed ball $B^{n}$, with $M (\infty)$
homeomorphic to $S^{n-1}$.  

\subsubsection{}\label{ss:equivalence}
Now the authors want to apply this idea to a locally symmetric space
$Y=\Gamma \backslash X$.  The geometry here is more intricate: there
are geodesic rays that don't go cleanly off to infinity, but close up
to form immersed loops.  Even worse, there are geodesics that reenter
a fixed compact set infinitely often.  This is already visible when $Y
= \Gamma \backslash \fH$ is a modular curve.  If $\gamma\subset \fH $
is a geodesic ray tending to an ideal point $\alpha\in \R$, then the
behavior of the image of $\gamma$ in $Y$ depends subtly on the
basepoint $\gamma (0)$ and the arithmetic nature of $\alpha$, in
particular the continued fraction expansion of $\alpha$.

Because of this phenomenon, the authors define a \emph{distance
minimizing ray} (DM ray) to be a geodesic ray $\gamma \colon \Rgeqz
\rightarrow Y$ that gives an isometric embedding of $\Rgeqz$ in $Y$.
More generally, they introduce \emph{eventually distance minimizing
geodesics} (EDM geodesics); by definition $\gamma \colon \R
\rightarrow Y$ is EDM if there exists $t_{0}>\!\!>0$ such that
$\gamma$ restricted to $\R_{\geq t_{0}}$ is DM.  These will be the
basic objects determining points on $Y (\infty)$.

Next they need an appropriate equivalence relation on geodesics.
Since parallelism makes no sense in $Y$, they generalize the
characterization of parallel rays in $\R^{n}$ through distances as in
\S\ref{ss:georn}.  They define two DM rays $\gamma ,\gamma '$ to be
\emph{equivalent} if
\[
\lim_{t\rightarrow \infty} \sup d (\gamma (t),\gamma ' (t))<\infty.
\]
In other words, as $t\rightarrow \infty$, the distances between the
corresponding points on $\gamma ,\gamma '$ remain bounded above. 

Finally they define $Y (\infty)$ as a set to be the set of DM rays
modulo equivalence.  After defining an appropriate topology on $Y\cup
Y (\infty)$, they obtain a compact Hausdorff space such that each DM
ray converges to the point on $Y (\infty)$ corresponding to its
equivalence class.  The resulting space is called the \emph{geodesic
compactification} of $Y$, and is denoted $Y^{\geo}$. 

\subsubsection{}
The basic idea of Ji--MacPherson's construction makes sense for any
complete noncompact Riemannian manifold, not just a locally symmetric
space, and in fact the geodesic compactification appears in the
literature in different guises.  For example, the geodesic
compactification of a Hadamard manifold is known as the \emph{conic
compactification} \cite{bgs}.(\footnote{A manifold is \emph{Hadamard}
if it is simply-connected, nonpositively curved, and complete
\cite{bgs}.  For example, $\fH\times \R$, $\fH \times \fH$, and
$\fH_{3}$ are all Hadamard.})

For complete noncompact Riemannian manifolds $M$, one has a very
general compactification technique due to Gromov \cite{bgs}: one
embeds $M$ in a space of continuous functions using the distance
function associated to the metric, and then takes the closure.  Ji and
MacPherson show that, if $M$ is a locally symmetric space, the
geodesic compactification coincides with Gromov's.

\subsection{More structure theory}
\subsubsection{}
Our next goal is a description of the geometry of $Y^{\geo}$, as well
as the relationship of $Y^{\geo}$ to $Y^{\borelserre}$ and
$Y^{\reductiveborelserre}$.  This requires substantially more notation
from the theory of algebraic groups
(\S\S\ref{ss:parabolic.start}--\ref{ss:parabolic.stop}).  We give
examples of most of the following in \S\ref{ss:slnexample} for $\bG =
\SL_{n}$; the inexpert reader may wish to skip ahead.

\subsubsection{}\label{ss:parabolic.start}
Let $\bS$ be a maximal $\Q$-split torus of $\bG$, with character
lattice (respectively, dual lattice of one-parameter subgroups) $X
(\bS)_{\Q}$ (resp., $X^{\vee}(\bS)_{\Q}$).  We denote by
$\scal{\phantom{a}}{\phantom{b}}\colon X (\bS)_{\Q}\times X^{\vee}
(\bS)_{\Q}\rightarrow \Z$ the natural unimodular pairing.

Let $\fg$ be the Lie algebra of $\bG$, and let $ \Phi = \Phi (\bG ,
\bS)$ be the roots of $\bS$ in the adjoint action on $\fg$.  Then
$\Phi$ is a root system in the vector space $X (\bS)_{\R} := X
(\bS)_{\Q}\otimes \R$.  This root system determines a hyperplane
arrangement $\{H_{\alpha}\mid \alpha \in \Phi \}$ in $X^{\vee}(\bS)_{\R}$ by
\[
H_{\alpha} = \{y\in X^{\vee}(\bS)_{\R}\mid  \scal{\alpha}{y} = 0\}.
\]
The connected components of the complement of this arrangement are
called \emph{Weyl chambers}.  The Weyl group $W = W (\bG )$ acts on $X^{\vee}
(\bS )_{\R}$ by transitively permuting the chambers.  If we fix a
chamber $C$, we determine a subset $\Phi^{+}\subset \Phi$ of positive
roots, and a subset $\Delta \subset \Phi^{+}$ of simple roots.

\subsubsection{}\label{ss:tits}
A closed subgroup $\bP \subset \bG$ is called \emph{parabolic} if it
contains a maximal connected solvable subgroup, and is called
\emph{$\Q$-parabolic} (or \emph{rational parabolic}) if it is defined
over $\Q$.

The \emph{spherical Tits building} $\sB  = \sB (\bG)$ is the simplicial complex
constructed as follows.  Simplices of $\sB$ are in bijection with proper rational
parabolic subgroups of $\bG$.  The vertices of $\sB$ correspond to the
maximal parabolic subgroups; given a set of such $\bP_{1},\dotsc
,\bP_{k}$, if the intersection $\bQ = \bP_{1}\cap \dotsb \cap \bP_{k}$ is a
rational parabolic subgroup, then there is a simplex in $\sB$ corresponding to
$\bQ$ with vertices corresponding to $\bP_{1},\dotsc ,\bP_{k}$.  It is
known that $\sB$ has pure dimension equal to $q-1$, where $q$ is the
$\Q$-rank of $\bG$.

The group of rational points $\bG (\Q)$ acts on $\sB$ through
conjugation of rational parabolic subgroups.  If $\Gamma \subset \bG
(\Q)$ is an arithmetic subgroup, then $\Gamma$ acts on $\sB$, and the
quotient $\Gamma \backslash \sB$ has simplices in bijection with the
$\Gamma$-conjugacy classes of rational parabolic subgroups.  It is
known that $\Gamma \backslash \sB$ is always finite.

\subsubsection{}\label{sss:levi}
Let $\bN_{\bP}\subset \bP$ be the unipotent radical of $\bP$.  The
quotient $\bL_{\bP} := \bP /\bN_{\bP}$ is called the \emph{Levi
quotient}.  It is a reductive group, defined over $\Q$ if $\bP$ is.

Let $\bS_{\bP}$ be a maximal $\Q$-split torus in the center of
$\bL_{\bP}$.  Let $A_{\bP} = S_{\bP}^{0}$, the connected component of
the identity in the group $S_{\bP}$ of real points 
$\bS_{\bP} (\R)$.  Put 
\[
\bM_{\bP} = \bigcap_{\alpha \in X (\bL_{\bP})_{\Q}} \ker \alpha^{2}.
\]
The group $\bM_{\bP}$ is semisimple, in general not connected.

There is a unique lift $i\colon \bL_{\bP}\rightarrow \bP$ compatible
with our choice of maximal compact subgroup $K\subset G$.  This
induces lifts of $\bS_{\bP}$ and $\bM_{\bP}$ to $\bP$ that allow us
to view the groups of real points $A_{\bP}$ and $M_{\bP}$ as subgroups
of $P$.  This leads to the \emph{Langlands
decomposition} of $P$:
\[
P = N_{\bP}A_{\bP}M_{\bP}.
\]
One can show that the map $N_{\bP}\times A_{\bP}\times
M_{\bP}\rightarrow P$ induced by multiplication is a diffeomorphism.

\subsubsection{}\label{ss:horo}
The Langlands decomposition leads to coordinates on the global
symmetric space $X = G/K$ as follows.  The group $P$ acts transitively on $X$.
Any $x\in X$ can be written as $x = uamx_{0}$, where $u\in
N_{\bP}$, $a\in A_{\bP }$, $m\in M_{\bP}$ are uniquely determined, and
where $x_{0}\in X$ is the basepoint determined by $K$.  Let $K_{\bP} =
M_{\bP}\cap K$ and put $X_{\bP} = M_{\bP}/K_{\bP}$.  Then $X_\bP$ is
the product of a global symmetric space of noncompact type with a
possible Euclidean factor.  The decomposition $x = uamx_{0}$ induces a
diffeomorphism 
\[
N_{\bP}\times X_{\bP}\times A_{\bP}\longrightarrow X
\]
by 
\[
(u,mK_{\bP}, a)\longmapsto uamx_{0}.
\]

\subsubsection{}\label{ss:parabolic.stop}
Let $\bP$ be a rational parabolic subgroup, and let $\fa_{\bP}$,
$\fn_{\bP}$ be the Lie algebras of $A_{\bP}$, $N_{\bP}$.  We denote by
$\Phi^{+} (\bP , A_{\bP})$ the positive roots of the adjoint action of
$\fa_{\bP}$ on $\fn_{\bP}$.  Put 
\begin{equation}\label{eq:aplusinfinity}
A_{\bP}^{+} (\infty) = \{H\in \fa_{\bP}\mid \alpha (H)>0, \quad 
(H,H)=1,  \quad \alpha \in \Phi^{+} (\bP , A_{\bP})\},
\end{equation}
where $(\phantom{a},\phantom{a})$ is the Killing form on $\fa_{\bP}$.
Let $\ol{A}_{\bP}^{+} (\infty)$ be the closure of $A_{\bP}^{+}
(\infty)$ obtained by replacing the conditions $\alpha (H)>0$ in
\eqref{eq:aplusinfinity} with $\alpha (H)\geq 0$.  Then
$\ol{A}_{\bP}^{+} (\infty)$ is homeomorphic to a closed simplex.
Define a simplicial complex $B (X)$ by 
\begin{equation}\label{eq:real1}
B (X) = \bigcup_{\bP} \ol{A}_{\bP}^{+} (\infty) / \sim,
\end{equation}
where the union is taken over all proper rational parabolic subgroups,
and we identify 
$\ol{A}_{\bP}^{+} (\infty)$ with a face of $\ol{A}_{\bQ}^{+} (\infty)$
if the two parabolic subgroups $\bP$, $\bQ$ satisfy $\bQ \subset
\bP$.  Then $B (X) \simeq \sB$, in other words the complex $B (X)$
is a realization of the Tits building $\sB$.  If $\Gamma$ is an
arithmetic subgroup, we similarly define 
\begin{equation}\label{eq:real2}
B (\Gamma \backslash X) = \bigcup_{\bP} \ol{A}_{\bP}^{+} (\infty) / \sim,
\end{equation} 
where the union is now taken over all $\Gamma$-conjugacy classes of
proper rational parabolic subgroups.  We have $B (\Gamma \backslash
X) \simeq \Gamma \backslash \sB$.

\subsubsection{}\label{ss:slnexample}
We consider the example $\bG = \SL_{n}$.

The maximal $\Q$-split torus $\bS$ is the subgroup of all diagonal
matrices.  The standard choice of positive roots is the set $\Phi^{+}
= \{e_{i}-e_{j}\mid 1\leq i<j\leq n \}$, where $\{e_{i} \}$ is the
standard basis of $\R^{n}$.  The simple roots $\Delta \subset
\Phi^{+}$ are the points $\{e_{i}-e_{i+1}\mid i=1,\dotsc ,n-1 \}$, and
the root lattice $X (\bS)_{\Q}$ can be identified with
$\{(x_{1},\dotsc ,x_{n})\in \Z^{n}\mid \sum x_{i} = 0 \}$.

Any proper $\Q$-parabolic subgroup is conjugate over $\SL_{n} (\Q)$ to
a standard proper $\Q$-parabolic subgroup.  The latter are indexed by
ordered positive partitions $\pi$ of $n$ with at least $2$ parts.
Given such a partition $\pi = (\pi_{1},\dotsc ,\pi_{k})$, the
corresponding parabolic subgroup has real points
\[
P = \Biggl\{\Biggl(\begin{array}{ccc}
P_{1}&\cdots&*\\
&\ddots&\vdots\\
0&&P_{k}
\end{array} \Biggr)\Biggm| P_{i}\in \GL _{\pi _{i}} (\R), \prod \det (P_{i}) = 1\Biggr\}. 
\]
The Langlands decomposition $P=N_{\bP}A_{\bP}M_{\bP}$ is given as
follows:
\begin{itemize}
\item $M_{\bP}\subset P$ is the subgroup of block diagonal matrices
such that each
block is an element of $\SL ^{\pm }_{\pi _{i}} (\R)$.  Here the $\pm$
means to take matrices with determinant $\pm 1$.
\item $A_{\bP}\subset P$ is the subgroup of block diagonal matrices
such that each $i$th block has the form $a_{i}I_{\pi _{i}}$, where
$a_{i}>0$ and $I_{\pi _{i}}$ is the $\pi _{i}\times \pi _{i}$ identity
matrix.
\item  $N_{\bP}\subset P$ is the subgroup such that each $i$th block equals 
$I_{\pi _{i}}$. 
\end{itemize}

The parabolic subgroup $\bP_{0}$ corresponding to the partition
$(1,\dotsc ,1)$ is called the Borel subgroup.  The Lie algebra
$\fa_{\bP_{0}}$ of $A_{\bP_{0}}$ can be identified with $X^{\vee}
(\bS)_{\R}$.  There are $n!$ Weyl chambers, each of which is an open
simplicial cone of dimension $n-1$.  The subset
$\overline{A}_{\bP_{0}}^{+} (\infty)$ is the intersection of the
closure of one of these cones with a sphere.  For any other standard
rational parabolic subgroup $\bQ$, we have an inclusion
$\fa_{\bQ}\subset \fa_{\bP_{0}}$ such that $\overline{A}_{\bQ}^{+}
(\infty)$ is identified with a proper face of
$\overline{A}_{\bP_{0}}^{+} (\infty)$. 

\begin{figure}[htb]
\psfrag{a}{$\alpha_{1}$}
\psfrag{b}{$\alpha_{2}$}
\psfrag{c}{$\alpha_{3}$}
\psfrag{Ainf}{$\overline{A}_{\bP_{0}}^{+} (\infty)$}
\psfrag{112}{$112$}
\psfrag{121}{$121$}
\psfrag{211}{$211$}
\psfrag{31}{$31$}
\psfrag{13}{$13$}
\psfrag{22}{$22$}
\begin{center}
\includegraphics[scale=0.5]{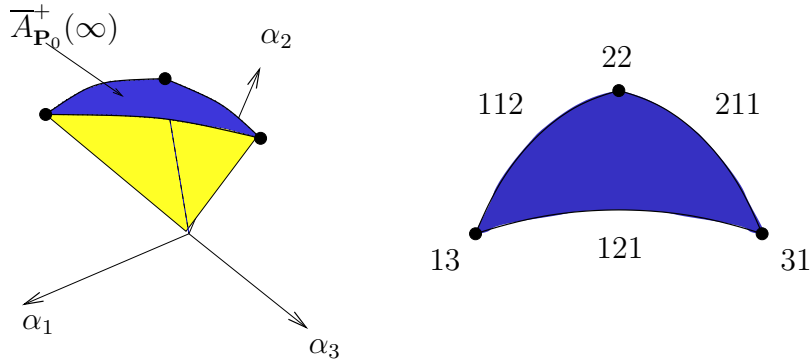}
\end{center}
\caption{The region $\ol{A}^{+}_{1111}= \overline{A}_{\bP_{0}}^{+}
(\infty)$ for $\SL_{4}$, where $\bP_{0}$ is a Borel subgroup.  \label{Ainf.fig}}
\end{figure}

Figure \ref{Ainf.fig} shows the situation for $\SL_{4}$.  The left
shows $\overline{A}_{\bP_{0}}^{+} (\infty)$ as the dark spherical
triangle topping off the simplicial cone.  In this figure we have
used the Killing form to identify $\fa_{\bP_{0}}$ with its dual, which
allows us to view the simple roots $\alpha_{i}$ alongside the cone.
Each simple root is orthogonal to a facet of the cone.
We have also rescaled the form in the definition
\eqref{eq:aplusinfinity} to make the picture clearer.

There are $7$ partitions $\pi$: $1111$, $211$, $121$, $112$, $31$,
$22$, $13$.  Here we abbreviate $(\pi_{1},\dotsc ,\pi_{k})$ by
$\pi_{1}\dotsb \pi_{k}$, and in what follows denote a parabolic
subgroup by its partition and drop $(\infty)$ from the notation.
The right of Figure \ref{Ainf.fig} shows how the faces of
$\ol{A}_{1111}^{+} = \ol{A}_{\bP_{0}}^{+} (\infty)$ are indexed by the partitions.  (To go
from the left figure to the right rotate the back of the dark
triangle forward.)  The edges of $\ol{A}^{+}$ correspond to partitions
with $3$ parts; the bottom is $\ol{A}_{121}^{+}$, the left
$\ol{A}_{112}^{+}$, and the right $\ol{A}_{211}^{+}$.  The vertices
correspond to partitions with $2$ parts; clockwise from the top they
are $\ol{A}_{22}^{+}$, $\ol{A}_{31}^{+}$, and $\ol{A}_{13}^{+}$.

\subsection{The geodesic compactification}
\subsubsection{}\label{ss:edmgeo}
We are now ready to discuss the geometry of $Y^{\geo}$.  Choose a
rational parabolic subgroup $\bP$, and write $X = N_{\bP}\times
X_{\bP}\times A_{\bP}$ as in \S\ref{ss:horo}.  If we choose $u\in
N_{\bP}$, $z\in X_{\bP}$, $a\in A_{\bP}$, and $H\in A_{\bP}^{+}
(\infty)$, we get a geodesic $\gamma$ on $X$ by 
\begin{equation}\label{eq:geodesic}
\gamma (t) = (u,z,a \exp (tH)) \subset N_{\bP}\times
X_{\bP}\times A_{\bP}, \quad t \in \R ,
\end{equation}
where $\exp\colon \fa_{\bP} \rightarrow A_{\bP}$ is the exponential
map.  The authors show that the image of $\gamma$ in $Y=\Gamma
\backslash X$ is an EDM geodesic, and that in fact \emph{any} EDM
geodesic on $Y$ has the form \eqref{eq:geodesic} for appropriately
chosen $\bP$.  Moreover, two geodesics $\gamma ,\gamma '$ of the form
\eqref{eq:geodesic} project to the same geodesic in $Y$, up to
reparametrization, if and only if $H=H'$, $\log a - \log a'$ is a
multiple of $H$, and $(u,z)= g(u',z')$ for some $g \in  \Gamma_{\bP}
:= \Gamma \cap \bP (\Q)$.   Hence all EDM geodesics on $Y$ have an especially
simple form. 

\subsubsection{}\label{sss:equivalence}
Next the authors investigate equivalence.  It turns out that the
equivalence class of the geodesic $\gamma (t) = (u,z,a \exp (tH))$
depends only on the $H$ component, and not $u$, $z$, or $a$.  This,
together with the explicit realizations
\[
B (X) \simeq \sB, \quad B (\Gamma \backslash
X) \simeq \Gamma \backslash \sB
\]
from \S\ref{ss:parabolic.stop}, shows that the boundary $\partial
Y^{\geo}$ is homeomorphic to $\Gamma \backslash \sB$.

This result motivates Ji and MacPherson to define another natural
compactification of $Y$: the \emph{Tits compactification} $Y^{\tits}$.
By definition, the (partial) Tits compactification of $X$ is the union
$X\cup \sB$, appropriately topologized so that the quotient $Y\cup
(\Gamma \backslash \sB)$ is compact and Hausdorff.  In $Y^{\tits}$,
the image of the geodesic $(u,z,a \exp (tH))$ converges to the point
in $\Gamma \backslash \sB$ corresponding to $H$ via the realization
\eqref{eq:real2}.

\subsection{The Borel--Serre and the reductive
Borel--Serre compactifications}
\subsubsection{}
The Tits compactification is a new compactification of $Y$ constructed
using the group theory of $\bG$, and as such properly joins the list
in \S\ref{ss:zoo}.  As we shall shortly see, $Y^{\tits}$ is as far
from $Y^{\borelserre}$ as possible: if the $\Q$-rank of $\bG$ is $>1$,
then the greatest common quotient of $Y^{\tits}$ and $Y^{\borelserre}$
is the one-point compactification of $Y$!  With hindsight, it is clear
that this must be the case, since the boundary of $Y^{\tits}\simeq
Y^{\geo}$ is constructed using rays that converge to each other at
infinity, and on $Y^{\borelserre}$ the metric degenerates.  Hence
equivalence classes of EDM geodesics cannot detect points in $\partial
Y^{\borelserre}$.

\subsubsection{}\label{ss:bsandrbs}
To understand the exact relationship between $Y^{\geo}$ and
$Y^{\borelserre}$, $Y^{\reductiveborelserre}$, we must first say more
about the construction of the latter spaces.  We begin with the
globally symmetric space $X$ and construct partial compactifications
$X^{\borelserre}$, $X^{\reductiveborelserre}$ by gluing on boundary
components for each proper rational parabolic subgroup.  Given such a
subgroup $\bP$, write $X = N_{\bP}\times X_{\bP}\times A_{\bP}$.  We
define boundary components $e^{\borelserre }(\bP)$,
$e^{\reductiveborelserre } (\bP )$ by
\[
e^{\borelserre } (\bP) = N_{\bP}\times X_{\bP}, \quad
e^{\reductiveborelserre } (\bP) = X_{\bP},
\]
and set 
\begin{equation}\label{eq:parcompact}
X^{\borelserre} = X \bigcup_{\bP} e^{\borelserre } (\bP), \quad
X^{\reductiveborelserre} = X \bigcup_{\bP} e^{\reductiveborelserre} (\bP),
\end{equation}
where the unions are taken over all proper rational parabolic
subgroups.

The topologies on \eqref{eq:parcompact} can be defined by specifying
the convergence properties of sequences
of points in $X$ and in the boundary components.  For instance, to
converge to a point on $e^{\borelserre } (\bP)$ starting from inside $X$, take a
sequence $\{y_{n} \} = \{(u_{n}, z_{n}, \exp (H_{n})) \} \subset N_{\bP}\times
X_{\bP}\times A_{\bP}$, and suppose that
\begin{enumerate}
\item $u_{n}\rightarrow u_{\infty}\in
N_{\bP}$, 
\item $z_{n}\rightarrow z_{\infty} \in X_{\bP}$, and 
\item for any $\alpha \in \Phi^{+} (\bP , A_{\bP})$, we have $\alpha
(H_{n})\rightarrow \infty$ as $n\rightarrow \infty$.
\end{enumerate}
Then in $X^{\borelserre}$, the sequence $y_{n}$ converges to
$(u_{\infty}, z_{\infty})\in e^{\borelserre } (\bP)$.  Similarly, in
the reductive Borel--Serre $X^{\reductiveborelserre}$, the same
sequence converges to the point $z_{\infty}\in e^{\reductiveborelserre
} (\bP)$.  Other sequences can be constructed that converge from one
boundary component to a point in another.  For a full list of all the
sequences that need to be specified to produce the correct topologies
on $X^{\borelserre}$ and $X^{\reductiveborelserre}$, we refer to
\cite[\S 7]{ji.macp}. 

Appropriately topologized, the quotients by the action of any
arithmetic group $\Gamma$ become the compact Hausdorff
compactifications $Y^{\borelserre}$ and $Y^{\reductiveborelserre}$.
We have a commutative diagram
\[
\cd{X^{\borelserre}}{Y^{\borelserre}}{X^{\reductiveborelserre}}{Y^{\reductiveborelserre}}.
\]
The left vertical map is the identity on $X$, and for each proper
$\Q$-parabolic subgroup $\bP$ is the projection $N_{\bP}\times
X_{\bP}\rightarrow X_{\bP}$.  The right vertical map is the identity
on $Y$; to understand it on the boundary components, we need more
notation.  Let $\Gamma_{N_{\bP}} = \Gamma \cap N_{\bP}$, and let
$U_{\bP}$ be the nilmanifold $\Gamma_{N_{\bP}}\backslash N_{\bP}$.
Let $\Gamma_{\bP} = \Gamma \cap P$, and let $\Gamma_{M_{\bP}} \subset
M_{\bP}$ be the discrete group obtained by projecting $\Gamma_{\bP}$
to $M_\bP$ via $P = N_{\bP}\times M_{\bP}\times A_{\bP} \rightarrow
M_{\bP}$.  Write $Z_{\bP}^{\borelserre}$ (respectively,
$Z_{\bP}^{\reductiveborelserre}$) for the boundary component of
$Y^{\borelserre}$ (resp., $Y^{\reductiveborelserre}$) corresponding to
$\bP$.  Then $Z^{\reductiveborelserre}_{\bP}$ is isomorphic to the
locally symmetric space $Y_{\bP} = \Gamma_{M_{\bP}}\backslash X_{\bP}$.  In
$Y^{\borelserre}$, the component $Z^{\borelserre}_\bP$ has the
structure of a fiber bundle over $Y^{\bP}$
with fiber $U_{\bP}$.  For each $\bP$, the map $Z^{\borelserre}_{\bP}\rightarrow
Z^{\reductiveborelserre}_{\bP}$ collapses the nilmanifold $U_{\bP}$ to
a point.

\subsubsection{}\label{ss:Nequivalence}
From the description of EDM geodesics on $Y$ given in
\S\ref{ss:edmgeo}, it is clear that all points in the boundaries $\partial
Y^{\borelserre}$ and $\partial Y^{\reductiveborelserre}$ can be
reached by following an EDM geodesic to its limit point.  Hence there
should be some construction of $Y^{\borelserre}$ and
$Y^{\reductiveborelserre}$ by putting a suitable equivalence relation
on the EDM geodesics.  This is indeed the case, and is carried out here.

If the $\Q$-rank of $\bG$ is $1$, the construction is quite simple.
For the Borel--Serre, Ji and MacPherson prove that there is a
bijection between the set of all EDM geodesics and points in $\partial
Y^{\borelserre}$, given by taking the endpoint of an EDM geodesic $\gamma$:
\[
\gamma \longmapsto \lim_{t\rightarrow \infty} \gamma (t) 
\]

The reductive Borel--Serre $Y^{\reductiveborelserre}$ is a quotient of
the Borel--Serre, so we need an equivalence relation on the EDM
geodesics.  We say two EDM geodesics $\gamma ,\gamma '$
are \emph{$N$-equivalent} (notation: $\gamma \Neq \gamma'$) if 
\[
\lim_{t\rightarrow \infty } d(\gamma (t), \gamma ' (t)) =0 
\]
for appropriate parametrizations of $\gamma , \gamma '$.  Denote the
$N$-equivalence class of $\gamma$ by $[\gamma]_{N}$.  Then
Ji--MacPherson show that there is bijection between the set of EDM
geodesics modulo $N$-equivalence and the points in $\partial
Y^{\reductiveborelserre}$, again given by taking the endpoint:
\begin{equation}\label{eq:rbsendpoint}
[\gamma]_{N} \longmapsto \lim_{t\rightarrow \infty} \gamma (t).
\end{equation}
Here the limit in \eqref{eq:rbsendpoint} is taken in
$Y^{\reductiveborelserre}$.

\subsubsection{}
For higher $\Q$-ranks, similar results hold, although we need more
equivalence relations to state them.  Let $\gamma$ be an EDM geodesic.
We define two sets $C (\gamma)$, $F (\gamma)$ of EDM geodesics related to $\gamma$ by
\begin{align*}
C (\gamma) &= \bigl\{\gamma '\bigm| \text{$d (\gamma (t), \gamma '(t))$ is
constant for $t>\!\!>0$} \bigr\}\\
F (\gamma) &= \bigl\{\gamma '\bigm|\lim_{t\rightarrow \infty} \sup d
(\gamma (t), \gamma ' (t))<\infty \bigr\}.
\end{align*}
The set $C (\gamma)$ is called the \emph{congruence bundle} of
$\gamma$.  It consists of all the EDM geodesics that are eventually at
constant distance to $\gamma$.  The set $F (\gamma)$ is called the
\emph{finite bundle} of $\gamma$.  Note that $F (\gamma)$ is actually
the set of all EDM geodesics that are \emph{equivalent} to $\gamma$ in
the sense used to construct the geodesic compactification
(\S\ref{ss:equivalence}). 

\subsubsection{}\label{ss:congbundle}
The congruence bundle $C (\gamma)$ can be turned into a metric space
with metric $\delta$ as follows.  First we put $\delta (\gamma
',\gamma) = c$, where $c$ is the constant in the definition of $C
(\gamma)$.  This constant can also be recovered as $\lim_{t\rightarrow
\infty} \bar d (\gamma (t), \gamma ')$, where $\bar d (\gamma (t),
\gamma ')$ is defined by
\[
\inf \bigl\{d (\gamma (t), \gamma ' (s))\bigm|s\in \R\bigr\}.  
\]
Via this description, we can also extend $\delta$ to all of $C
(\gamma)$.  The authors show that $(C (\gamma), \delta)$ is complete,
and in fact has the following concrete form.  If 
\[
\gamma (t) =
(u,z,a\exp (tH))\subset N_{\bP}\times
X_{\bP}\times A_{\bP},
\]
then
\begin{equation}\label{eq:Cgamma}
C (\gamma) \simeq Y_{\bP}\times \Span (H)^{\perp}.
\end{equation}
Here $Y_{\bP}$ is as in \S\ref{ss:bsandrbs},
and $\Span (H)^{\perp}$ is the orthogonal complement to the line
through $H$ in $\fa_{\bP}$.

\subsubsection{}
Now we are ready to define our next equivalence relation.  Let
$\gamma$ be EDM.  We define the \emph{rank} $r (\gamma)$ of $\gamma$
by 
\[
r=r (\gamma) = \bigl\{k\in \Z\bigm| \text{there exists a faithful
isometric action of $\R^{k-1}$ on $C (\gamma)$}\bigr\}.
\]
We then say $\gamma '\in C (\gamma)$ is \emph{$L$-related} to $\gamma$
(notation: $\gamma \Leq \gamma '$) if $\gamma$, $\gamma '$ belong to
the same $\R^{r-1}$-orbit.  Here the $L$ stands for
\emph{linear}; one pictures the $\R^{r-1}$-action as
linearly sliding the geodesics in $C (\gamma)$ around in the $\Span
(H)^{\perp}$ factor from \eqref{eq:Cgamma}.  

The authors show that $L$-equivalence extends to an equivalence
relation on the set of all EDM geodesics.  In terms of
\eqref{eq:Cgamma}, $L$-equivalence can be written as follows.  If
$\gamma_{1} \Leq \gamma_{2}$, then there exists $\gamma $ with
$\gamma_{i}\in C (\gamma)$.  Write $C (\gamma)$ as in
\eqref{eq:Cgamma}, and write $\gamma_{i} (t) = (u_{i}, z_{i},
a_{i}\exp (tH_{i}))$.  Then $\gamma_{1}\Leq \gamma_{2}$ means that the
$u_{i}$ project to the same point in the nilmanifold $U_{\bP}$ and the
$z_{i}$ to the same point in the locally symmetric space
$Y_\bP$.  From this it also follows that $r
(\gamma)$ is the $\Q$-rank of $\bP$.

\subsubsection{}
The restriction of $L$-equivalence to the finite bundle $F (\gamma)$
also induces an equivalence relation.  The dimension of the quotient
$F (\gamma)/\Leq$ is called the \emph{mobility degree} of $\gamma$,
and is denoted $\mu (\gamma)$. 

Using the mobility degree we can define another equivalence relation
on EDM geodesics, called \emph{$R$-equivalence}.   The $R$ stands for
\emph{rotation}.  Let $\gamma_{0}, \gamma_{1}$ be EDM and let
$[\gamma_{i}]_{L}$ be their $L$-equivalence classes.  We say
$[\gamma_{0}]$ is $R$-equivalent to $[\gamma_{1}]$, and write
$[\gamma_{0}]_{L}\Req [\gamma_{1}]_{L}$, if there exists a family
$\gamma_{s}$, $s\in [0,1]$ of
EDM geodesics interpolating $\gamma_{0}$ and $\gamma_{1}$ with the
following properties:
\begin{enumerate}
\item $d (\gamma_{s_{1}} (t), \gamma_{s_{2}} (t)) = c |s_{1}-s_{2}|t$
for $t\geq 0$, for some constant $c$ and for all $s_{1},s_{2}\in
[0,1]$, and 
\item the mobility degree $\mu$ restricted to $\gamma_{s}$ is constant.
\end{enumerate}
Two EDM geodesics $\gamma,\gamma '$ are said to be \emph{$RL$-related}
(notation: $\gamma \RLeq \gamma '$) if their $L$-classes are
$R$-related: $[\gamma]_{L}\Req [\gamma']_{L}$. 

For the reductive Borel--Serre, we will need to combine
$RL$-equivalence with $N$-equivalence from \S \ref{ss:Nequivalence}.
We say $\gamma$ and $\gamma '$ are \emph{$NRL$-equivalent} if there
exists an EDM geodesic $\gamma ''$ such that $\gamma \RLeq \gamma ''$
and $\gamma ' \Neq \gamma ''$.

\subsubsection{}
We can finally explain how to construct $Y^{\borelserre}$ and
$Y^{\reductiveborelserre}$ using EDM geodesics.  If the $\Q$-rank of
$\bG$ is bigger than $1$, then $\partial Y^{\borelserre}$ is in
bijection with the $RL$-equivalence classes of EDM geodesics, via the
endpoint map:
\[
[\gamma]_{RL} \longmapsto \lim_{t\rightarrow \infty} \gamma (t),
\]
where the limit is taken in $Y^{\borelserre}$.  
For the reductive Borel--Serre, the boundary $\partial
Y^{\reductiveborelserre}$ is in bijection with the $NRL$-equivalence
classes of EDM geodesics by the
endpoint map:
\[
[\gamma]_{NRL} \longmapsto \lim_{t\rightarrow \infty} \gamma (t),
\]
where the limit is now taken in $Y^{\reductiveborelserre}$.

\subsubsection{}\label{ss:geotobs}
We can now finally explain the relationship between $Y^{\geo}$ and the
compactifications $Y^{\borelserre}$, $Y^{\reductiveborelserre }$.
First we need an explicit realization of $RL$-equivalence. 

Suppose $\gamma (t) = (u,z,a\exp (tH))$ and $ \gamma \RLeq \gamma '$.
Then it turns out we can write $\gamma ' (t)$ in the form 
\[
\gamma '(t) = ( u,z,a'\exp
(tH'))
\]
for the same parabolic subgroup $\bP$.  Hence two geodesics are
$RL$-equivalent if their $N_{\bP}$ and $X_{\bP}$ components coincide;
the $A_{\bP}$ part is irrelevant.

Notice that this is exactly the \emph{opposite} of the basic
equivalence relation used to construct $Y^{\geo}$: for the geodesic
compactification, the $A_{\bP}$ component is the \emph{only} component
of an EDM geodesic that matters; the $N_{\bP}$ and $X_\bP$ components
play no role (\S\ref{sss:equivalence}).  This is the sense in which
$Y^{\geo}$ is as far as possible from $Y^{\borelserre}$.

\subsubsection{}
We conclude our discussion by giving fanciful pictures of the local
structure of compactifications $Y^{\borelserre}$,
$Y^{\reductiveborelserre}$, and $Y^{\geo}$ for $\bG =\SL_{3}$ near the
boundary components corresponding to the $\Gamma$-conjugacy classes of
the standard rational parabolic subgroups.

\begin{figure}[htb]
\psfrag{xg}{$Y$}
\psfrag{xp0}{$Y_{111}$}
\psfrag{xp1}{$Y_{{12}}$}
\psfrag{xp2}{$Y_{{21}}$}
\psfrag{n1}{$U_{12}$}
\psfrag{n2}{$U_{21}$}
\psfrag{n12}{$U_{111}$}
\begin{center}
\includegraphics[scale=0.2]{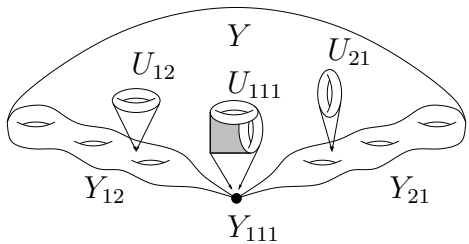}
\end{center}
\caption{Borel--Serre compactification for $\SL_{3}$\label{bs.fig}}
\end{figure}

Figure \ref{bs.fig} shows the Borel--Serre compactification.  This
figure is based on one by MacPherson \cite{mac.course}.  We use the
notation of \S\ref{ss:slnexample}, so that partitions are $111$, $21$,
and $12$.  If $\pi$ is a partition corresponding to the standard
rational parabolic subgroup $\bP$ with Langlands decomposition $P =
N_{\bP}A_{\bP}M_{\bP }$, we write $Y_{\pi}$ for the locally symmetric
space $\Gamma_{M_{\bP}}\backslash X_\bP $ and $U_{\pi}$ for the
nilmanifold $\Gamma_{N_{\bP}}\backslash N_{\bP}$
(\S\ref{ss:bsandrbs}).  The codimension 1 corners are both torus
bundles over locally symmetric spaces for $\SL_{2}$.  In these cases
the unipotent radicals $N_{12}$, $N_{21}$ are both isomorphic to
$\R^{2}$ with the subgroups $\Gamma_{12}$, $\Gamma_{21}$ isomorphic to
$\Z ^{2}$ acting by translation.  The alignment of the tori in the
figure indicates the structure of the unipotent radicals.  For
example, the unipotent radical $N_{12}$ consists of all real matrices
of the form
\[
\left(\begin{array}{ccc}
1&*&*\\
0&1&0\\
0&0&1
\end{array} \right).
\]

In the codimension 2 corner $Z_{111}$ the locally symmetric space
$Y_{111}$ is a point.  Hence all the topology of $Z_{111}$ is contained in the
nilmanifold $U_{111}$.  This 3-manifold is known as the
\emph{Heisenberg manifold}.  It can be written as a bundle
$T^{2}\rightarrow Z_{111} \rightarrow S^{1} $ in two different ways,
reflecting the two subgroups $N_{12}, N_{21}\subset
N_{111}$.(\footnote{Incidentally, both of these bundles are different
from the $T^{2}$-bundle in \S\ref{ss:quadraticexample}.})

\begin{figure}[htb]
\psfrag{xg}{$Y$}
\psfrag{xp0}{$Y_{111}$}
\psfrag{xp1}{$Y_{12}$}
\psfrag{xp2}{$Y_{21}$}
\begin{center}
\includegraphics[scale=0.2]{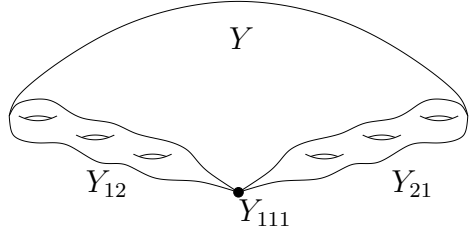}
\end{center}
\caption{Reductive Borel--Serre compactification for $\SL_{3}$\label{rbs.fig}}
\end{figure}

Figure \ref{rbs.fig} shows the reductive Borel--Serre
compactification.  In this figure the nilmanifold fibers have been
collapsed to points, and all that remains in the boundary are the
lower-rank locally symmetric spaces.

\begin{figure}[htb]
\psfrag{xg}{$Y$}
\psfrag{xp0}{$\sB_{111}$}
\psfrag{xp1}{$\sB_{12}$}
\psfrag{xp2}{$\sB_{21}$}
\begin{center}
\includegraphics[scale=0.2]{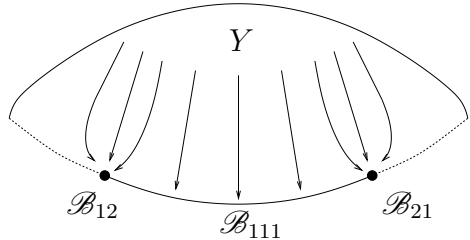}
\end{center}
\caption{Geodesic compactification for $\SL_{3}$\label{g.fig}}
\end{figure}

Finally, Figure \ref{g.fig} shows the geodesic compactification with
some geodesics suggestively converging to boundary components.  Here
we write $\sB_{\pi}$ for the simplex in the building corresponding to
the standard parabolic subgroup $P_{\pi}$.  As expected from
\S\ref{ss:geotobs}, the boundary components here are in some sense
totally opposite to those of $Y^{\borelserre}$ and
$Y^{\reductiveborelserre}$: as the codimension increases in the
boundaries of $Y^{\borelserre}$ and $Y^{\reductiveborelserre}$, so
does the \emph{dimension} of the corresponding boundary components of
$Y^{\geo}$.  This illustrates why, if the $\Q$-rank is $>1$, the
greatest common quotient of $Y^{\borelserre}$ and $Y^{\geo}$ is the
one-point compactification.

\subsection{Complements}

\subsubsection{}
We conclude by briefly summarizing some of the other material in
\cite{ji.macp}.

\subsubsection{}\label{sss:laplace}
Let $M$ be a complete Riemannian manifold, and let $\Delta$ be the Laplace
operator on $L^{2} (M)$.  It is known that if $M$ is noncompact, then
the continuous spectrum of $\Delta$ (if it exists) cannot change under
compact perturbation of $M$.  Hence one has the natural problem of
trying to understand the connection between the continuous spectrum of
$\Delta$ and compactifications of $M$.

In the case that $M$ is a locally symmetric space $Y=\Gamma \backslash
X$, it is known that $\Delta$ has continuous spectrum thanks to
Langlands's study of Eisenstein series \cite{langlands, mw}.  Ji and
MacPherson are able to describe the continuous spectrum of $\Delta$ in
terms of the geometry of $Y^{\geo}$ by reinterpreting Langlands's
fundamental work.  This provides a very accessible and geometric
introduction to the intricate constructions in \cite{langlands, mw}.

\subsubsection{}
 Let $M$ be a complete Riemannian manifold and let $\lambda$ be any
real number less than the bottom of the spectrum of the Laplace
operator $\Delta$.  Associated to this data is a certain
compactification of $M$, the \emph{Martin compactification}
$M^{\martin}_{\lambda}$.  The precise definition is somewhat involved;
we refer to \cite{ji.etal} and \cite[\S 15]{ji.macp} for details.
Using Eisenstein series, Ji and MacPherson show that for a locally
symmetric space $Y$ the geodesic compactification $Y^{\geo}$
canonically injects into the Martin compactification
$Y^{\martin}_{\lambda }$.  They also conjecture that this injection is
in fact a homeomorphism.

\bibliographystyle{amsplain_initials}
\bibliography{bobmacp}

\end{document}